\documentclass[12pt]{amsart}\usepackage{amssymb,verbatim,pstricks}
\newcommand{\nocom}{}
\nocom
\newtheorem{thm}{Theorem}[section] 
\newtheorem{pro}{Proposition}[section]
\newtheorem{lem}{Lemma}[section]
\newtheorem{cor}{Corollary}[section]

\theoremstyle{definition}
\newtheorem{dfn}{Definition}[section]

\theoremstyle{remark}
\newtheorem{rem}{Remark}

\begin{document}
\title[Hecke algebras and $3$-manifolds]
{Hecke algebras, modular categories and $3$-manifolds quantum invariants}
\author{Christian Blanchet}
\address{Université de Nantes-UMR 6629, 2 rue de la Houssinière,
BP92208, 44322 NANTES Cedex 3, FRANCE}
\email{blanchet$\char'100$math.univ-nantes.fr}
\date{Febrary 1998}

\begin{abstract}
We construct modular categories from
Hecke algebras at roots of unity.
For a special choice of the framing parameter, we recover
the Reshetikhin-Turaev invariants of closed $3$-manifolds
constructed from the quantum groups $U_q sl(N)$ by
 Reshetikhin-Turaev and Turaev-Wenzl, and from skein theory
 by Yokota.
 The possibility of such a construction was suggested by Turaev,
as a consequence of Schur-Weil duality.
 We then discuss the choice of the framing parameter. This leads,
for any rank $N$ and  level $K$, to a modular category
$\widetilde{\mathrm{H}}^{N,K}$ and a reduced  invariant 
$\tilde\tau_{N,K}$. If $N$ and $K$ are coprime, then this invariant
coincides with the known invariant $\tau^{PSU(N)}$ at level $K$.
If $\mathrm{gcd}(N,K)=d>1$, then we show that the reduced invariant admits
spin or cohomological refinements, with a nice decomposition formula
which extends a theorem of H. Murakami.
\end{abstract}
\maketitle
\section*{Introduction}
Our goal is to construct modular categories from the
 Hecke algebras of type $A$ at roots of unity, and to obtain
certain reductions and refinements of them.
Our main results are the followings.
\begin{itemize}
\item We give  a (reasonably self contained)
construction of  modular categories underlying
the known $SU(N)$ quantum invariants.
As expected, at rank-level $(N,K)$, isomorphism classes of simple objects
are indexed by the set of Young diagrams with at most $N-1$ lines and
$K$ columns, whose cardinality is $\frac{(N+K-1)!}{(N-1)!K!}$.
 The dimensions of the TQFT modules are given by the Verlinde
 formula.
\item At rank-level $(N,K)$ we obtain a reduced modular category
$\widetilde{\mathrm{H}}^{N,K}$,
and a reduced invariant $\tilde \tau_{N,K}$. Here the number of
non-isomorphic simple objects is
$\frac{(N+K-1)!}{N!K!}d$, with $d=\mathrm{gcd}(N,K)$.
Except for $\mathrm{gcd}(N,K)=1$, in which case we recover the $PSU(N)$
invariant,
 this result seems to be new.
\item We use a grading modulo $d=\mathrm{gcd}(N,K)$ 
to define invariants of $3$-manifolds equipped
with what we call a  spin structure
with coefficients modulo $d$, and we prove a decomposition
formula which extends a theorem of H. Murakami \cite{Mu}.
This refinement holds when $d$ is even, and $\frac{N}{d}$,
$\frac{K}{d}$ are odd (the spin case);
in the other cases, we obtain cohomological refinements.
\end{itemize}

The quantum invariants,
predicted by Witten using Chern-Simons theory and path integrals,
were first constructed by Reshetikhin-Turaev \cite{RT} and Turaev-Wenzl
\cite{TW1} using representation theory of quantum groups.
 The work of Turaev \cite{Turaev} shows that a key concept
 in the construction of these invariants, as well as in extending them
 to a  topological quantum field theory (TQFT), is
 that of a modular category.
 
 A modular category is a braided category with some additional
algebraic features (duality, twist, a finite set of simple objects
satisfying a domination property and a non-degeneracy axiom).
The interest of this concept is that it provides
a topological quantum field theory in dimension $3$, and
in particular, invariants of links and $3$-manifolds. However,
in particular examples, it is not  easy to
define precisely the modular category and
to check the required properties.
  For the $SU(N)$ invariants, first constructed using representation
 theory of $U_q(sl(N))$, it is known that underlying modular categories
can be derived from the category of representations of
the quantum group at roots of unity \cite{TW1,MW,An}.
  In section \ref{sectionSUN}, we will give an alternative elementary
  construction of modular   categories producing the same
invariants\footnote
{Two modular categories producing the same invariant should be equivalent.
  A complete treatment of this question has yet to be formed.}.

A skein theoretic construction 
of the $SU(N)$ invariants was obtained by Yokota \cite{Yo};
subsequent developments towards the associated TQFT's were
made by Lickorish \cite{LiSamp,LiTQFT}.
Our main tool here will be this skein method
combined with the structure of the Hecke algebra,
which in our context is obtained as the Homflypt skein
module of a cylinder $\mathbf{D}^2\times[0,1]$ with
boundary points.
 This algebra
 has been intensively studied (see \cite{Geck} for a list of
 references). For our purpose, we emphasize the work of
 Jones \cite{Jones}
and Wenzl \cite{Wenzl};
see also \cite{Gyo}, \cite{DJ}.  Our  normalization and description of
idempotents coincide
 with  those of Aiston-Morton
 \cite{Aiston,AM,MA}.

Following the work of Kirby and Melvin for the $SU(2)$
case \cite{KM}, Kohno and Takata have studied symmetry formulas
for the $SU(N)$ quantum invariants, and defined the $PSU(N)$
invariants \cite{KT1,KT2}.
 This was used by H. Murakami in \cite{Mu}.
 Our reductions and refinements formulas generalize
 these results.
 
 Invariants of $3$-manifolds from  Hecke algebras
 were obtained by Wenzl in \cite{Wenzl2};
 modular and semi-simple categories from unoriented link invariants
(BCD case) are considered by Turaev and Wenzl in \cite{TW2}.

Masbaum and Wenzl \cite{MW} have proved  the integrality of $SU(N)$
quantum invariants at roots of unity of prime order, and
they have shown that this follows from existence of integral
modular categories. This is developped by Brugui\`eres in \cite{Brug}.

\noindent{\bf Acknowledgments.} The author wishes to thanks 
A. Beliakova, A. Brugui\`eres, N. Habegger
G. Masbaum,  V. Turaev  and P. Vogel
for useful discussions and suggestions.

\section{Young idempotents and Homflypt skein theory}\label{homflyCAT}
\subsection{The Homflypt functor}
Let $M$ be an oriented $3$-manifold.
We denote by $\mathcal{H}(M)$  the
$k$-module freely generated by
isotopy classes of framed links
in $M$ with the Homflypt relations 
given in figure \ref{homflyrel}.
 Here a framing  is a trivialization of the normal
 bundle; this gives an orientation of the link.
\begin{figure}[h]
\begin{center}
\vspace{5pt}
$a^{-1}$
\begin{pspicture}[.4](0,0)(1,1)
\psline{->}(.1,0)(.9,1)
\psline{-}(.9,0)(.6,.4)
\psline{->}(.4,.6)(.1,1)
\end{pspicture}
$\ -\ a$
\begin{pspicture}[.4](0,0)(1,1)
\psline{->}(.9,0)(.1,1)
\psline{-}(.1,0)(.4,.4)
\psline{->}(.6,.6)(.9,1)
\end{pspicture}
$\ =\ (s-s^{-1})$
\begin{pspicture}[.4](0,0)(1,1)
\pscurve{->}(.1,0)(.4,.5)(.1,1)
\pscurve{->}(.9,0)(.6,.5)(.9,1)
\end{pspicture}\\[10pt]
\begin{pspicture}[.4](0,0)(1,1)
\pscurve{-}(.1,0)(.15,.45)(.45,.75)(.7,.5)(.4,.2)(.20,.38)
\pscurve{->}(.14,.62)(.1,.9)(.1,1)
\end{pspicture}
$\ =\ av^{-1}$
\begin{pspicture}[.4](0,0)(1,1)
\psline{->}(.5,0)(.5,1)
\end{pspicture}
\hspace{2cm}
\begin{pspicture}[.4](0,0)(1,1)
\pscurve{<-}(.1,1)(.15,.55)(.45,.25)(.7,.5)(.4,.8)(.20,.62)
\pscurve{-}(.14,.38)(.1,.1)(.1,0)
\end{pspicture}
$\ =\ a^{-1}v$
\begin{pspicture}[.4](0,0)(1,1)
\psline{->}(.5,0)(.5,1)
\end{pspicture}\\[10pt]
$L\ \cup\ $
\begin{pspicture}[.4](0,0)(1,1)
\pscircle(.5,.5){.4}
\end{pspicture}
$\ =\ {\frac{v^{-1}-v}{s-s^{-1}}}\ L$
\vspace{5pt}
\end{center}
\caption{\label{homflyrel} Homflypt relations}
\end{figure}

Here $k$ is an integral domain containing the invertible elements
$a$, $v$, $s$;
we suppose moreover that ${s-s^{-1}}$ is invertible in $k$.
 (For $L\neq \emptyset$, the third equality is a consequence of
the others.)

The Homflypt polynomial \cite{homfly,PT} gives an isomorphism
$$\begin{array}{lrcl}
\langle \dots \rangle :&\mathcal{H}({\bf S}^3)&\rightarrow &k\\
&L&\mapsto &\langle L\rangle\end{array}$$
normalized by $\langle \emptyset\rangle=1$.

An embedding $j: M\hookrightarrow N$ gives a well defined operator
$\mathcal{H}(j): \mathcal{H}(M)\rightarrow \mathcal{H}(N)$. 
 This makes $\mathcal{H}$ into a functor from the category
 of $3$-manifolds whose morphisms are isotopy classes of embeddings, to
 the category of $k$-modules.

 To an oriented embedding of a disjoint union of solid tori
$$g=\coprod_{i=1}^mg_i : \ \coprod_{i=1}^m {\bf D}^2_i
\times {\bf S}_i^1\rightarrow {\bf S}^3$$
is associated a multilinear map
$${\mathcal{H} }(g):
\ {\mathcal{H}}({\bf D}^2\times {\bf S}^1)^{\otimes m}
\rightarrow {\mathcal H}({\bf S}^3)\approx k$$
This map only depends  on the isotopy class of the framed link 
$L=(L_1,\dots,L_m)$ underlying $g$.
The image of $x_1\otimes\dots\otimes x_m$ under this map
 is said to be obtained by cabling the components
 $L_i$ with the skein elements $x_i$,
and is denoted by
$\langle L_1(x_1),\dots,L_m(x_m)\rangle$ or
$\langle L(x_1,\dots,x_m)\rangle$.

In the case where $M$ has non-empty boundary, we may fix a framed link
$l$ in the boundary of $M$, and define the relative
skein module $\mathcal{H}(M,l)$.
By comparing with the orientation of $\partial M$
(our convention for orientation of the boundary is
outgoing normal vector first), the framing of $l$ gives
an orientation of its components.

\subsection{The Hecke category}
The {\em $k$-linear Hecke category} ${\mathrm{H}}$ is defined as follows.
An object in this category is a disc $\mathbf{D}^2$ equipped with a framed
link
(note that isotopy is not allowed).
 If $\alpha=(\mathbf{D}^2,l_0)$ and $\beta=(\mathbf{D}^2,l_1)$
are two objects,
 the module $Hom_{\mathrm{H}}(\alpha,\beta)$ is
$\mathcal{H}( \mathbf{D}^2\times[0;1],l_0\times 0\amalg l_1\times 1)$.
 The notation ${\mathrm{H}}(\alpha,\beta)$ and ${\mathrm{H}}_\alpha$ will be used respectively
for $Hom_{\mathrm{H}}(\alpha,\beta)$ and $End_{\mathrm{H}}(\alpha)$.
For composition, we use the {\em covariant} notation
$$\begin{array}{rcl}
{\mathrm{H}}(\alpha,\beta)\times {\mathrm{H}}(\beta,\gamma)&\longrightarrow &{\mathrm{H}}(\alpha,\gamma)\\
(f,g)&\longmapsto &fg
\end{array}$$
When we draw a figure to describe a morphism, the {\it time} parameter
goes upwards, so that the morphism $fg$ is depicted
with $g$ lying above $f$.

\noindent{\bf Note.} The terminology {\em Hecke category} was introduced
by Turaev in \cite{T1}. His definition gives a category
equivalent to an unframed version of ours.

We will use the notation $\hat{f}$ for the closure
in $\mathcal{H}(\mathbf{D}^2\times\mathbf{S}^1)$
of the morphism $f\in {\mathrm{H}}_\alpha$. Let $U_0$ be a $0$-framed unknot
in $\mathbf{S}^3$. We use the notation
$\langle f \rangle$ for $\langle U_0(\hat f) \rangle$
(the quantum trace).

 We will simply denote  by $n$
the object formed by $n$ standard 
 points (along the real axis)
 equipped with the standard
framing $(1,\sqrt{-1})$.

 By using a standard embedding
${\bf D}^2\amalg {\bf D}^2\hookrightarrow {\bf D}^2$, we make
 ${\mathrm{H}}$ into a  monoidal
category\footnote{The associativity
isomorphisms just move the embedded discs along the real axis,
and will be omitted.  Note that the object
$1^{\otimes n}$ is defined up to these associativity isomorphisms,
and is canonically isomorphic to the standard object $n$.}.
 This category has a braiding and a twist operator. We can
define in ${\mathrm{H}}$ a duality rule so that we get a ribbon category
(\cite{Turaev}, see also \cite{Vo}). We proceed as follows.
To an object $\alpha=(\mathbf{D}^2,l)$, we associate
the object $\alpha^*=(\mathbf{D}^2,-\overline l)$, where
$-\overline l$ is the framed link obtained by applying to $l$
the differential of $z\mapsto -\overline z$, and define the morphisms
$b_\alpha\in {\mathrm{H}}(0,\alpha\otimes \alpha^*)$
and $d_\alpha\in {\mathrm{H}}(\alpha^*\otimes\alpha,0)$ according to
the figure below (a copy of $1_\alpha$ is embedded
along the framed arc).
\centerline{
        \begin{pspicture}(2.5,2)
\put(-0.8,.9){$b_\alpha=$}
\pscurve{->}(1.5,1.5)(1.47,1.3)(1,.8)(.53,1.3)(.5,1.5)
\put(0.15,1.2){$\alpha$}
\end{pspicture}
\hspace{1cm}
\begin{pspicture}(2.5,2)
\put(-0.2,.9){$d_\alpha=$}
\pscurve{->}(2,.5)(1.97,.7)(1.5,1.2)(1.03,.7)(1,.5)
\put(2.1,.5){$\alpha$}
\end{pspicture}}

Our purpose is to discuss which modular categories
arise from this ribbon Hecke category.
 For this, we need 
  a finite set of  simple objects with nice properties
  (\cite{Turaev} p74).

\subsection{Idempotents in the Hecke algebra}
 The algebra ${\mathrm{H}}_n$
 is isomorphic to the quotient of the algebra
of the braid group $k[B_n]$ by the Homflypt relation

\begin{center}
\vspace{5pt}
$a^{-1}$
\begin{pspicture}[.4](0,0)(1,1)
\psline{->}(.1,0)(.9,1)
\psline{-}(.9,0)(.6,.4)
\psline{->}(.4,.6)(.1,1)
\end{pspicture}
$\ -\ a$
\begin{pspicture}[.4](0,0)(1,1)
\psline{->}(.9,0)(.1,1)
\psline{-}(.1,0)(.4,.4)
\psline{->}(.6,.6)(.9,1)
\end{pspicture}
$\ =\ (s-s^{-1})$
\begin{pspicture}[.4](0,0)(1,1)
\pscurve{->}(.1,0)(.4,.5)(.1,1)
\pscurve{->}(.9,0)(.6,.5)(.9,1)
\end{pspicture}
\vspace{5pt}
\end{center}
which is the Hecke algebra of type $A_{n-1}$.
 An unframed version of this result was proved independently by
 Morton-Traczyk \cite{MT} and Turaev \cite{T1}.
This algebra is known to be generically
semi-simple.
 It is a deformation of the algebra of the symmetric group,
 and its structure can be obtained by extending
 the classical  Young theory
\cite{Bou,Jones,Wenzl,Gyo,DJ}.

 Recall that, to a partition of $n$,
 $\lambda=(\lambda_1\geq\dots\geq\lambda_p\geq1)$,
 $\lambda_1+\dots+\lambda_p=n$, is associated a Young diagram
 of size $|\lambda|=n$, which we denote also by $\lambda$. This diagram
 has $n$ cells indexed by
$\{(i,j),\ 1\leq i \leq p,\ 1\leq j\leq \lambda_i\}$.

If $c$ is the cell of index $(i,j)$ in a Young diagram,
its hook-length $hl(c)$ and its content $cn(c)$
are defined by
$$hl(c)=\lambda_i+\lambda_j^\vee-i-j+1,\ \ cn(c)=j-i.$$
Here $\lambda^\vee$ is the transposed Young diagram and
$\lambda^\vee_j$ is the length of the $j$-th column
of $\lambda$ (the $j$-th line of $\lambda^\vee$).

For $n\geq0$, the quantum integer $[n]$, and the
quantum factorial $[n]!$ are defined by
$[n]=\frac{s^{n}-s^{-n}}{s-s^{-1}}$ and $[n]!=\prod_{j=1}^n [j]$.\\
For a Young diagram $\lambda$,
we will use the notation
$[hl(\lambda)]$, for the product over all cells
of the quantum
hook-lengths.
$$[hl(\lambda)]=\prod_{\mathrm{\scriptstyle cells}} [hl(c)]$$

\noindent{\bf{Symmetrizers.}}
Let $\sigma_i\in {\mathrm{H}}_n$, $i=1,\dots,n-1$, be represented by
the standard generators of the braid group $B_n$
(the strand numbered $i$ crosses over the strand numbered
$i+1$).
\begin{pro}
If $[n]!$ is invertible in $k$, then there exists a unique idempotent
$f_n\in {\mathrm{H}}_n$ such that $\forall i\ \sigma_if_n=as f_n=f_n\sigma_i$,
 and a unique idempotent
$g_n\in {\mathrm{H}}_n$ such that $\forall i\ \sigma_ig_n=-as^{-1} g_n=g_n\sigma_i$.
\end{pro}
\begin{proof}
It is shown in \cite{Morton,AM}  that the deformation $f_n$
of the Young symmetrizer, given below, satisfies the required condition.
$$f_n=\frac{1}{[n]!}s^{-\frac{n(n-1)}{2}}
\sum_{\pi\in \mathcal{S}_n} (as^{-1})^{-l(\pi)}w_\pi$$
Here $w_\pi$ is the positive braid associated with the
permutation $\pi$, and $l(\pi)$ is the length of $\pi$.

One can also construct $f_n$ recursively using the formulas below (\cite{Yo}).
$$[2]f_2=
s^{-1}\begin{pspicture}[.4](0,0)(1,1)
\pscurve{->}(.1,0)(.4,.5)(.1,1)
\pscurve{->}(.9,0)(.6,.5)(.9,1)
\end{pspicture}
+a^{-1}
\begin{pspicture}[.4](0,0)(1,1)
\psline{->}(.1,0)(.9,1)
\psline{-}(.9,0)(.6,.4)
\psline{->}(.4,.6)(.1,1)
\end{pspicture}$$
$$[n+1]f_{n+1}=-[n-1]f_n\otimes 1_1 +[2][n]
(f_n\otimes 1_1)(1_{n-1}\otimes f_2)(f_n\otimes 1_1)$$
\noindent(Here, $1_p\in {\mathrm{H}}_p$ is the identity.)

Suppose now that
$f'_n$ satisfies the condition of the theorem. Then
 we have that $f_nf'_n=f'_n=\zeta f_n$, and $\zeta=1$
  by idempotence.

 For $g_n$ we can proceed similarly,
 either with the deformed antisymmetrizer
 $$g_n=\frac{1}{[n]!}s^{\frac{n(n-1)}{2}}
\sum_{\pi\in \mathcal{S}_n} (-as)^{-l(\pi)}w_\pi$$
 or with the recursive formulas
 $$g_1=1_1$$
 $$g_{n+1}=1_1\otimes g_n -\frac{[2][n]}{[n+1]}
(1_1\otimes g_n)( f_2\otimes 1_{n-1})(1_1\otimes g_n)$$
\end{proof}
Note that $1_2=f_2+g_2$; this can be used to obtain
more symmetric recursive formulas for $f_n$ and $g_n$.
Our choice  immediately gives the following  lemma \cite{Yo},
which is useful in Yokota's skein computations.
\begin{lem}
For any two integers $p$, $q$ such that $[p+1]!$ and $[q+1]!$
are invertible, one has
\begin{eqnarray*}
[p+q]f_p\otimes g_q&=&[p+1][q]
(1_p\otimes g_q)( f_{p+1}\otimes 1_{q-1})(1_p\otimes g_q)\\&&
+[p][q+1](f_p\otimes 1_q)(1_{p-1}\otimes g_{q+1})(f_p\otimes 1_q)
\end{eqnarray*}
\end{lem}
\noindent{\bf Aiston-Morton description of Young symmetrizers.}
 For a Young diagram $\lambda$ of size $n$, we denote
by $\square_\lambda$ the object of the category ${\mathrm{H}}$ formed
with one point  for each cell $c$ of $\lambda$; if $c$ has index $(i,j)$
($i$-th line, and $j$-th column),
then the corresponding point in $\mathbf{D}^2$ is
$\frac{j+i\sqrt{-1}}{n}$.
Suppose that $\lambda=(\lambda_1\geq\dots\geq\lambda_p\geq1)$, and
that $\lambda^\vee=(\lambda_1^\vee\geq\dots\geq\lambda_q^\vee\geq1)$
is the transposed Young diagram.\\
Let $F_\lambda$ (resp. $G_\lambda$) be the element
in $ {\mathrm{H}}_{\square_{\lambda}}$ formed with one copy
of $[\lambda_i]!f_{\lambda_i}$ along the line $i$, for $i=1,\dots,p$
(resp. one copy
of $[\lambda^\vee_j]! g_{\lambda_j^\vee}$ along the column $j$,
 for $j=1,\dots,q$).
 Note that these expressions have no denominators.

 In the proposition below, $<$ denotes the lexicographic
 ordering, and $\lambda$ and $\mu$ are Young diagrams with
 the same number of cells.
\begin{pro}\label{minimal}
a) If $\mu<\lambda$, then
$F_\lambda  {\mathrm{H}}(\square_{\lambda},\square_{\mu})
 G_\mu=0$.\\
b) If $\lambda<\mu$, then
$G_\lambda  {\mathrm{H}}(\square_{\lambda},\square_{\mu})
 F_\mu=0$.\\ 
 c) One has $F_\lambda  {\mathrm{H}}_{\square_{\lambda}}  G_\lambda=
kF_\lambda G_\lambda$.
\end{pro}
This proposition is proved in \cite{AM}, section 4.
We outline here the proof given there.
\begin{proof}
a) Every family of $|\lambda|!$ braids in
${\mathrm{H}}(\square_{\lambda},\square_\mu)$,
which induces all  bijections between  the cells of $\lambda$
and the cells of $\mu$, is a basis of
${\mathrm{H}}_{\square_{\lambda}}$.
 The hypothesis $\mu<\lambda$ implies that
any bijection between the cells of $\lambda$
and the cells of $\mu$ carries at least two cells in some line
of $\lambda$,  to cells in the same column of $\mu$.
 So we can find a basis of ${\mathrm{H}}(\square_{\lambda},\square_{\mu})$
 represented by braids which connect in a separate cylinder
  two points in some line of $\square_{\lambda}$ with two points
 in some column of $\square_{\mu}$.
 One can deduce that
$F_\lambda  {\mathrm{H}}(\square_{\lambda},\square_{\mu})
 G_\mu=0$. Assertion b) is shown similarly.
 
 c) We say that a permutation $\pi$ between the cells of $\lambda$
 does not separate, if some pair of cells in the same line
 is mapped  to some pair of cells in the same
 column; other permutations are said to separate.
Using that any separating permutation $\pi$ of 
the cells can be written $\pi=\pi_R\pi_C$ (permutations act on
the left), where $\pi_R$ (resp. $\pi_C$)
preserves the rows (resp. the columns), we see
that we can find a basis represented by braids $b_\pi$
indexed by permutations
such that
$$\begin{cases}
F_\lambda b_\pi G_\mu=0 &\text{ if $\pi$ does not separate,}\\
F_\lambda b_\pi G_\mu=\zeta_\pi F_\lambda  G_\mu&
\text{ if  $\pi$  separates.}
\end{cases}$$
Statement c) follows.
\end{proof}
The argument used in the proof of a) above shows the
following.
\begin{lem}
a) If $p>\lambda_1$, and for some object $\alpha$,
$x=x_1\otimes f_p\otimes x_2\in {\mathrm{H}}_\alpha$, then
$x{\mathrm{H}}(\alpha,\square_{\lambda})G_\lambda=0$.\\
b)  If $p>\lambda^\vee_1$, and for some object $\alpha$,
$x=x_1\otimes g_p\otimes x_2\in {\mathrm{H}}_\alpha$, then
$F_\lambda {\mathrm{H}}(\square_\lambda,\alpha)x=0$.
\end{lem}
Let $\tilde y_\lambda=F_\lambda G_\lambda$. 
 A consequence of  proposition \ref{minimal} is that
$\tilde y_\lambda$ is a quasi-idempotent.
\begin{pro}\label{quasi}
One has $\tilde y_\lambda^2=[hl(\lambda)]\tilde y_\lambda$.
 \end{pro}
 This proposition is proved by Yokota in \cite{Yo},
lemma 2.3 (see also \cite{Aiston2});
  the fact that our parameters $a$ and $v$ are more generic is irrelevant
  in this skein computation.\\
If $[hl(\lambda)]$ is invertible, we define the
 idempotent $y_\lambda$ by
$$y_\lambda=[hl(\lambda)]^{-1}\tilde y_\lambda$$
 Suppose that $\lambda$ and $\mu$ are Young diagrams with the
same number of cells,
such that $[hl(\lambda)]$ and $[hl(\mu)]$ are invertible.
From proposition \ref{minimal}, we have.
\begin{pro}
a) If $\mu\neq\lambda$, then
$y_\lambda  {\mathrm{H}}(\square_{\lambda},\square_{\mu})
 y_\mu=0$.\\
 b) One has that $y_\lambda  {\mathrm{H}}_{\square_{\lambda}}  y_\lambda=
ky_\lambda$.
\end{pro}
Let $\lambda\subset \mu$ be two Young diagrams,
the complement
of $\lambda$ in $\mu$ is called 
a skew Young diagram and is denoted by $\mu/\lambda$.
As above we can define the object
$\square_{\mu/\lambda}$ in the category ${\mathrm{H}}$ with one point for
each cell. The following is  proven in the same way
as proposition \ref{quasi}.
\begin{lem}\label{absorbe}
Let $\lambda\subset \mu$ be two Young diagrams,
 and let
$v\in {\mathrm{H}}_{\square_{\mu/\lambda}}$. One has that
 $$F_\mu(\tilde y_{\lambda}\otimes v)G_\mu=
[hl(\lambda)]F_\mu(1_{\square_\lambda}\otimes v)G_\mu $$
\end{lem}
We will need the following formulas.
\begin{pro}\label{braiding}
a) Suppose that $\lambda\subset\mu$ are two Young diagrams such that
the skew diagram $\mu/\lambda$ has only one cell $c$
($|\mu|=|\lambda|+1$). Then

\centerline{
        \begin{pspicture}(0,0)(8,4)
        \psframe(2,0.5)(3,1)
        \psline(2.5,0.25)(2.5,0.5)
        \psframe(2,3)(3,3.5)
        \psline{->}(2.5,3.5)(2.5,3.75)
        \put(2.2,0.7){$y_\mu$}
        \put(2.2,3.2){$y_\mu$}
        \psframe(1.5,1.5)(2.5,2)
        \put(1.7,1.7){$y_\lambda$}
        \psline[linearc=.3](2.2,1)(2.2,1.1)(2,1.4)(2,1.5)
        \psline[linearc=.3](2,2)(2,2.1)(2.1,2.6)
        \psline[linearc=.3](2.15,2.8)(2.2,2.9)(2.2,3)
        \pscurve(2.8,1)(2.8,1.5)(2.7,2)(2.2,2.2)
        \pscurve(1.9,2.3)(1.7,2.5)(2.8,2.9)(2.8,3)
\put(3.5,2){ $=a^{2|\lambda|}s^{2cn(c)} y_\mu$}
\end{pspicture}}
b) (Framing coefficient)
 
\centerline{
        \begin{pspicture}(0,0)(8,3)
        \psframe(1,0.5)(2,1)
        \psline(1.5,0.25)(1.5,0.5)
        \put(1.2,0.7){$y_\lambda$}
        \pscurve(1.5,1)(1.5,1.5)(1.65,1.7)(1.8,1.5)(1.7,1.2)(1.65,1.2)
        \pscurve{->}(1.35,1.2)(1.3,1.2)(1.15,1.5)(1.4,2)(1.5,2.5)
        \put(3,1.3)
{ $=a^{|\lambda|^2}v^{-|\lambda|}s^{2\sum_{cells}cn(c)}y_\lambda$}
\end{pspicture}}
\end{pro}
Statements a) and b) can be deduced from each other.
 Statement b) is theorem 5.5 of \cite{AM}; the proof given
there is an elementary
skein calculation which contains statement a). These results
can also be deduced from \cite{Wenzl2}. Statement b) is a framed
version of lemma 3.2.1 there.

\subsection{Structure of the generic Hecke algebra and path idempotents}
Here we suppose that $n$ is fixed, and that, in the domain $k$,
$[j]$ is invertible
for every $j\leq n$, so that all the idempotents
$y_\lambda$ exist in ${\mathrm{H}}_n$.

 A standard tableau $t$ with shape a Young diagram
 $\lambda=\lambda(t)$ is a labelling of the cells, with the integers
 $1$ to $n$, which is increasing along lines and columns.
 We denote by $t'$ the tableau obtained by removing
 the cell numbered by $n$.
 We define $\alpha_t\in {\mathrm{H}}(n,\square_\lambda)$ and
$\beta_t\in {\mathrm{H}}(\square_\lambda,n)$
 by
 $$\alpha_1=\beta_1=1_1$$
 $$\alpha_t=(\alpha_{t'}\otimes 1_1)\varrho_t y_\lambda$$
 $$\beta_t=y_\lambda \varrho_t^{-1}(\beta_{t'}\otimes1_1) $$
 Here $\varrho_t\in {\mathrm{H}}(\square_{\lambda(t')}\otimes 1,
\square_\lambda)$ is a standard  isomorphism which is obtained by
moving the added point, first along
a line parallel to the imaginary axis,
and then along a line parallel to the real axis.

Note that $\beta_\tau\alpha_t= 0$ if 
$\tau\neq t$, and
$\beta_t\alpha_t=y_{\lambda(t)}$.
\begin{pro}
The family $\alpha_t\beta_\tau$ for all standard tableaux $t,\tau$ such that
$\lambda(t)=\lambda(\tau)$ forms a basis for ${\mathrm{H}}_n$.
\end{pro}
\begin{proof}
We have (here $\delta$ is the Kronecker delta).
$${\alpha_t\beta_\tau\alpha_s\beta_\sigma}
=\delta_{\tau s} {\alpha_t\beta_\sigma}$$
This shows independence. Moreover the number of vectors
 is equal to the dimension.
  This shows the proposition over the field of quotients
  of the domain $k$.
   The above formula gives the coordinate forms and shows
   that the result is valid over $k$.
 \end{proof}
This gives explicitly the semi-simple decomposition of ${\mathrm{H}}_n$.
 The simple components are indexed by Young diagrams.
 The component indexed by $\lambda$ is the two-sided ideal
 generated by $y_\lambda$; its rank is $d_\lambda^2$, where
 $d_\lambda$ is the number of standard tableaux with shape $\lambda$.
The $\alpha_t\beta_\tau$ are {\em matrix units} in the sense of Ram-Wenzl
 \cite{RW}. The diagonal elements $p_t=\alpha_t\beta_t$
are the path idempotents described in
 \cite{Wenzl}; the minimal central idempotents are
 the $z_\lambda=\sum_{\lambda(t)=\lambda}p_t$.
 
\begin{thm}[Branching formula]\label{branch}
$$y_\lambda\otimes 1_1=
\sum_{\genfrac{}{}{0pt}{2}{\lambda\subset \mu}{|\mu|=|\lambda|+1}}
(y_\lambda\otimes 1_1)y_\mu (y_\lambda\otimes 1_1)$$
\end{thm}
\noindent We have omitted in this formula the standard
isomorphism between $\lambda\otimes 1$ and $\mu$.

This branching formula for the path idempotents is given
by
Wenzl in \cite{Wenzl}. A similar formula is given
by Yokota \cite{Yo}, prop. 2.11, with a proof using skein
calculus.
\begin{proof}
We continue to denote by $z_\mu$ the minimal central idempotent
corresponding to the simple component indexed by $\mu$
in ${\mathrm{H}}_{\square_\lambda\otimes 1}$.
 We have $$y_\lambda\otimes 1_1=\sum_\mu z_\mu(y_\lambda\otimes 1_1)
=\sum_{\genfrac{}{}{0pt}{2}{\lambda\subset \mu}{|\mu|=|\lambda|+1}}
(y_\lambda\otimes 1_1)z_\mu(y_\lambda\otimes 1_1)$$
From the formula $z_\mu=\sum_{\lambda(t)=\mu} p_t$,
we see that in the above, only those tableaux $t$
with $\lambda(t')=\lambda$ contribute,
and all these contributions are proportional
to $(y_\lambda\otimes 1_1)y_\mu (y_\lambda\otimes 1_1)$, by lemma
\ref{absorbe}. By using idempotence,
we have that $(y_\lambda\otimes 1_1)z_\mu(y_\lambda\otimes 1_1)=
(y_\lambda\otimes 1_1)y_\mu (y_\lambda\otimes 1_1)$.
\end{proof}
\begin{pro}[Quantum dimension]
$$ \langle y_\lambda\rangle=
\prod_{\mathrm{cells}}{\frac{v^{-1}s^{cn(c)}-vs^{-cn(c)}}
{s^{hl(c)}-s^{-hl(c)}}}$$
\end{pro}
Recall that $\langle y_\lambda\rangle$ is the Homflypt polynomial
of a $0$-framed unknot cabled with the closure $\hat y_\lambda$ of
$y_\lambda$ in $\mathbf{D}^2\times
\mathbf{S}^1$.\\
We will denote $\langle  y_\lambda\rangle$
simply by $\langle\lambda\rangle$,
and call it the quantum dimension of $\lambda$.

The assertion above can  be proven by  a skein calculation
(see proposition 2.4 in \cite{Yo} or \cite{Aiston2}). An alternative proof
(\cite{Wenzl}) is to use
the Young algebra.  From the branching formula  (\ref{branch}),
we can deduce that $\hat y_\lambda$ is the $\lambda$-indexed
Schur polynomial
in the $\hat g_k$; hence the general formula
follows from $\langle  g_k\rangle=\prod_{j=1}^k
\frac{v^{-1}s^{j-1}-vs^{1-j}}{s^j-s^{-j}}$, which can be proven by using the
recursive formula for $g_k$.

\subsection{The $\mathcal{C}$-completed Hecke category}
Suppose $\mathcal{C}$ is a set of Young diagrams $\lambda$, such
that the Young idempotents $y_\lambda$ exist.
The $\mathcal{C}$-completed Hecke category ${\mathrm{H}}^{\mathcal{C}}$
is defined as follows.

 An object in this category is a disc $\mathbf{D}^2$
 equipped with a labelled framed link $l=((l_1,\lambda^{(1)}),\dots,
 (l_m,\lambda^{(m)}))$, where $\lambda^{(1)},\dots,\lambda^{(m)}$
 are Young diagrams in $\mathcal{C}$.

If $\alpha=(\mathbf{D}^2,l)$ is such an object, its expansion
$E(\alpha)=(\mathbf{D}^2,E(l))$ is obtained by embedding
the object  $\square_\lambda$
in a neighbourhood of $l_i$, according to the framing. 
 The tensor product $y_{\lambda^{(1)}}\otimes\dots
\otimes y_{\lambda^{(m)}}$
defines an idempotent $\pi_\alpha\in {\mathrm{H}}_{E(\alpha)}$.
 The module ${\mathrm{H}}^{\mathcal{C}}(\alpha,\beta)$ is defined
 by $${\mathrm{H}}^{\mathcal{C}}(\alpha,\beta)=\pi_\alpha
{\mathrm{H}}(E(\alpha),E(\beta))\pi_\beta$$
 The duality extends to the category ${\mathrm{H}}^{\mathcal{C}}$ in a
natural way, and we again have  a ribbon category.
 We  denote simply by $\lambda$
the object of
${\mathrm{H}}^{\mathcal{C}}$ which is a disc $\mathbf{D}^2$ with the origin labelled
 by $\lambda$.

\subsection{ Homflypt calculus using ribbon graphs}
Following Turaev (\cite{Turaev}, I.2), we can define
the category $Rib_{{\mathrm{H}}}$ of ribbon graphs over ${\mathrm{H}}$, and use
the canonical functor $F_{\mathrm{H}}: Rib_{{\mathrm{H}}}\rightarrow {\mathrm{H}}$.
 A colored ribbon graph gives a morphism in the category ${\mathrm{H}}$
 (the functor $F_{\mathrm{H}}$ will be implicit).
  This can be described as follows. For each band in the graph,
  colored with $\alpha$,
  embed (using the framing of the band) a copy of $1_\alpha$,
the identity of $\alpha$;
  for each loop colored with  $\alpha$ embed a copy of
  $\hat1_\alpha$, the closure of $1_\alpha$;
  for each coupon colored with the morphism $f$, embed a copy of $f$.

 We can proceed similarly with the $\mathcal{C}$-completed
Hecke category
 ${\mathrm{H}}^\mathcal{C}$.
 
\section{The modular categories $\mathrm{H}^{SU(N,K)}$
and $\mathrm{H}^{PSU(N,K)}$}\label{sectionSUN}
\subsection{Roots of unity}
In this section,
we suppose that $s$ is a primitive $2(N+K)$-th root of unity,
and that $v=s^{-N}$ (rank $N$ and level $K$,
$N\geq 2$ and $K\geq1$).
We suppose moreover that $N+K$ is invertible in $k$.
A consequence is  that
$[n]$ is invertible for $n<N+K$. This can be seen as follows.
 We have that
$$[n]=s^{-n+1}\prod_{1<j|n} \phi_j(s^2)$$
 Here $\phi_j\in \mathbb{Z}[X]$ is the $j$-indexed cyclotomic polynomial.
The required invertibility is a consequence of the following 
lemma. Here we suppose $p\geq2$.
 \begin{lem}
 If $j\notin p\mathbb{Z}$, then $\phi_j$ divides $p$
in $\mathbb{Z}[X]/\phi_p$.
 \end{lem}
 \begin{proof}
 Let $d=\mathrm{gcd}(p,j)$. In $\mathbb{Z}[X]$, one has the
  relation $U(X^j-1)+V(X^p-1)=X^d-1$.
If $j>d$, this implies the relation
$U_1\phi_j+V_1\phi_p=1$, hence $\phi_j$ is invertible
in $\mathbb{Z}[X]/\phi_p$.
If $j=d$ then
using the  derivative of $X^p-1=\phi_jT$,
  we get $p\equiv XT'\phi_j$ mod $\phi_p$.
\end{proof}

We observe that
   the idempotent
$y_\lambda$ exists for every Young diagram
$\lambda$ with $\lambda_1+\lambda^\vee_1\leq N+K$.
 We denote by ${\mathcal C}_{N,K}$ the set of these Young diagrams
and consider the category $\mathrm{H}^{(N,K,a)}$
 obtained from the
${\mathcal C}_{N,K}$-completed Hecke category,
${\mathrm{H}}^{\mathcal{C}_{N,K}}$
 by applying
the following purification procedure (\cite{Turaev} p504).
\begin{dfn}
A morphism $f\in {\mathrm{H}}^{\mathcal{C}_{N,K}}(\alpha,\beta)$
is negligible iff
$$\forall g\in {\mathrm{H}}^{\mathcal{C}_{N,K}}(\beta,\alpha),
\ \ \langle {fg}\rangle=0$$
\end{dfn}

Objects of $\mathrm{H}^{(N,K,a)}$
  are those of ${\mathrm{H}}^{\mathcal{C}_{N,K}}$.
The module  $\mathrm{H}^{(N,K,a)}(\alpha,\beta)$ is the quotient
of ${\mathrm{H}}^{\mathcal{C}_{N,K}}(\alpha,\beta)$ by the submodule
of negligible morphisms.
 One can verify that composition and tensor product
are well defined on the quotient.
The negligible morphisms give local relations
in Homflypt modules. We denote by
$\mathcal{H}^{(N,K,a)}(M)$, the quotient of 
$\mathcal{H}^{\mathcal{C}_{N,K}}(M)$ by these relations.
 Note that for $M=\mathbf{S}^3$, no new relation appears.
In the case $M=\mathbf{D}^2\times\mathbf{S}^1$, the
algebra structure is well defined on the quotient.
 The multilinear form $\langle L(\dots)\rangle$
is well defined on
$\mathcal{H}^{(N,K,a)}
(\mathbf{D}^2\times\mathbf{S}^1)$.
 More generally, ribbon graphs may be colored using the 
 category $\mathrm{H}^{(N,K,a)}$.

For a Young diagram $\lambda=(\lambda_1,\dots,\lambda_p)
\in \mathcal{C}_{N,K}$, the identity morphism
$1_\lambda$ is negligible if and only if
its quantum dimension $\langle\lambda\rangle$ is zero.
 This will be the case if and only if
 $K<\lambda_1$ or
$N<\lambda_1^\vee$.
We will use the following sets of Young diagrams.
 $$\overline\Gamma_{N,K}=\{(\lambda_1,\dots,\lambda_p),\
\lambda_1\leq K\text{ and }p\leq N\}$$
 $$\Gamma_{N,K}=\{(\lambda_1,\dots,\lambda_p),\
\lambda_1\leq K\text{ and }p<N\}$$
Denote by $1^N$ (resp. $(K)$)
the diagram,
with one column containing $N$ cells
(resp. one line containing $K$ cells).
The  proposition below shows that
$1^N$ (resp. $(K)$) is practically 
ininfluential.
 It will also explains our choice of the framing parameter $a$.
\begin{pro}
The following identities hold in the category
$\mathrm{H}^{(N,K,a)}$

\centerline{
        \begin{pspicture}(0,0.5)(5,3.5)
        \psframe(1.5,1.5)(2.5,2)
        \put(1.7,1.7){$g_N$}
        \psline[linearc=.3](2,1)(2,1.5)
        \psline[linearc=.3](2,2)(2,2.6)
        \psline[linearc=.3]{->}(2,2.8)(2,3.2)
        \pscurve(2.8,1)(2.8,1.5)(2.7,2)(2.2,2.2)
        \pscurve{->}(1.9,2.3)(1.7,2.5)(2.75,2.9)(2.8,3.2)
\put(3.5,1.5){ $=\ \ a^{2N}s^2 $}
\end{pspicture}
        \begin{pspicture}(0,0.5)(5,3.5)
        \psframe(1.5,1.5)(2.5,2)
        \put(1.7,1.7){$g_N$}
        \psline[linearc=.3](2,1)(2,1.5)
        \psline[linearc=.3]{->}(2,2)(2,3.2)
        \psline{->}(2.8,1)(2.8,3.2)
        \end{pspicture}}
\centerline{
\begin{pspicture}(0,0)(5,4)
        \psframe(1.5,1.5)(2.5,2)
        \put(1.7,1.65){$f_K$}
        \psline[linearc=.3](2,1)(2,1.5)
        \psline[linearc=.3](2,2)(2,2.6)
        \psline[linearc=.3]{->}(2,2.8)(2,3.2)
        \pscurve(2.8,1)(2.8,1.5)(2.7,2)(2.2,2.2)
        \pscurve{->}(1.9,2.3)(1.7,2.5)(2.75,2.9)(2.8,3.2)
\put(3.5,1.5){ $=\ \ a^{2K}s^{-2} $}
\end{pspicture}
        \begin{pspicture}(0,0)(5,4)
        \psframe(1.5,1.5)(2.5,2)
        \put(1.7,1.65){$f_K$}
        \psline[linearc=.3](2,1)(2,1.5)
        \psline[linearc=.3]{->}(2,2)(2,3)
        \psline{->}(2.8,1)(2.8,3)
\end{pspicture}
}
\end{pro}
\begin{proof}
We justify the first equality.  The proof of the second one is similar.
The branching formula \ref{branch} gives ($g_{N+1}$ is negligible)
$$g_N\otimes 1_1=(g_N\otimes 1_1)y_{(2,1^{N-1})}(g_N\otimes 1_1)$$
The formula comes from  proposition \ref{braiding}.
\end{proof}
Observe that (at rank $N$, level $K$), we have
$\langle  g_N\rangle=1$, and that the framing coefficient
for $g_N$ is $(a^Ns)^N$.  Using this we obtain the following corollary.


\begin{cor}\label{column}
Suppose that $a^Ns=1$.  Then the morphism represented by
a colored ribbon graph is not changed
\begin{verse}
if some band colored with
the object $1^N$
is twisted or moved across any other band,\\
or if some loop colored with
the object $1^N$
is removed.
\end{verse}
\end{cor}

\noindent{\bf Until the end of this section, we suppose that}
 $\mathbf{a^Ns=1}$.
 If $\lambda$ is a Young diagram in
$\overline{\Gamma}_{N,K}$, we denote by
 $\lambda^\star$ the skew diagram $\mu/\lambda$,
 with $\mu=\lambda_1^N$; up to a rotation in the plane,
$\lambda^\star$ is a Young diagram in $\Gamma_{N,K}$.
\begin{lem}\label{qdim}
a) For any $\lambda\in \overline{\Gamma}_{N,K}$, one has
$\langle {\lambda^\star}\rangle=
\langle {\lambda}\rangle$.\\
b) For any $\lambda\in \overline{\Gamma}_{N,K}$ of the form
 $\lambda=1^N+\nu$, one has
$\langle L( \hat y_{\lambda},\dots)\rangle=
\langle L( \hat y_{\nu},\dots\rangle$.\\
In particular, for $j\leq K$, one has $\langle j^N\rangle=1$
\end{lem}
\begin{proof}
If $v=s^{-N}$ (rank $N$), the quantum dimension
is given by $$\displaystyle \langle \lambda \rangle=
\prod_{\mathrm{cells}}\frac{[N+cn(c)]}{[hl(c)]}.$$
\end{proof}
\begin{lem}\label{rect}
For any $\lambda\in \overline\Gamma_{N,K}$ , we have\\
 \centerline{
        \begin{pspicture}(0,-.5)(5,3.5)
        \psframe(1,1)(3,2)
        \psline{->}(1.5,2)(1.5,3)\psline(1.5,0)(1.5,1)
        \put(1.7,1.3){$\lambda_1^N$}
        \pscurve{->}(2.5,2)(2.6,2.3)(3,2.5)(3.5,1.5)
        \pscurve(3.5,1.5)(3,.5)(2.6,.7)(2.5,1)
        \put(3.5,1.9){$\lambda^\star$}
        \put(4,1){\large $=\frac{1}{\langle \lambda\rangle}y_\lambda$.}
\end{pspicture}
}
\end{lem}
\begin{proof}
From the definition of the idempotents,
we see that the left hand side is proportional to $y_\lambda$.
 The coefficient is obtained by comparing the quantum traces.
 \end{proof}
\noindent {\bf Reversing  orientation.}
Let $L$ be a framed link, and let $L'$ be the framed link
obtained from
$L$ by reversing the orientation of the first component
(i.e. by changing the sign of the second vector in the trivialization of
the normal bundle).
\begin{pro}\label{reverse}
For any $\lambda\in \overline{\Gamma}_{N,K}$, one has
$$\langle L(\hat y_\lambda,\dots)\rangle=
\langle L'(\hat y_{\lambda^\star},\dots)\rangle$$
\end{pro}
\begin{proof}
Using  lemma \ref{rect},
we introduce somewhere the idempotent associated with the diagram
$\lambda_1^N$. Then we can cut the band colored with 
$\lambda_1^N$, move one of the ends along the first component $L_1$
of $L$ and glue it back again, so that the band
colored with $\lambda^\star$ goes along $L_1$ with the reverse
orientation. By  corollary \ref{column}  this does not change the
evaluated Homflypt polynomial.
 Using  lemma \ref{rect} again, we get the result.
 \end{proof}
The above suggests that we could build an isomorphism
between the trivial object $0$ and $1^N$, and also
between the dual of $\lambda$ (denoted by
$\lambda^*$) and $\lambda^\star$
(note the difference between $\ast$ and $\star$). In order to do that,
we will add some morphisms.

We saw that we can define morphisms in the Hecke category,
and more generally skein elements in Homflypt skein
modules by using colored ribbon graphs.
We  extend the $\mathcal{C}_{N,K}$-completed
 Hecke category and the Homflypt skein theory
by allowing  
 incoming or outgoing vertices, colored with $r=1^N$.
Together with the relation of isotopy rel. boundary, the Homflypt relations,
and the negligible morphisms, we add the relation given
 by gluing an incoming $r$-colored vertex with an outgoing one.
(Note that the half twist is not trivial; hence the orientations
of the glued ends of the ribbon edges must be respected.)

We will denote by $\overline{\mathrm{H}}^{(N,K,a)}$
the category whose objects are those of $\mathrm{H}^{(N,K,a)}$
and whose morphisms are defined using the extended Homflypt skein
theory above, and by $\overline{\mathcal{ H}}^{(N,K,a)}$
the extended Homflypt skein functor.

From corollary \ref{column}, we deduce that
the Homflypt polynomial extends to an isomorphism
$\overline{\mathcal{H}}^{(N,K,a)}(\mathbf{S}^3)\approx k$.

In the category  $\overline{\mathrm{H}}^{(N,K,a)}$,
an object $\lambda\in \Gamma_{N,K}$ is still simple,
and its dual
 $\lambda^*$ is isomorphic to $\lambda^\star$.
 
From the branching formula (theorem \ref{branch}),
we can see that the simple objects in $\Gamma_{N,K}$ dominate
the category $\overline{\mathrm{H}}^{(N,K,a)}$.
This gives us all the defining properties of a
modular category except the non-degeneracy axiom.
We say that
$\overline{\mathrm{H}}^{(N,K,a)}$ is a {\em pre-modular category}.

Note that distinct objects $\lambda,\mu\in \Gamma_{N,K}$
are not isomorphic, but this does not imply
the non-degeneracy axiom.

\subsection{The handle slide condition}
We say that $\Omega\in{\mathcal{ H}}^{(N,K,a)}
({\bf D}^2\times {\bf S}^1)$ satisfies
the Kirby condition if
$$({\rm K})\ 
\begin{cases}
\forall x\in {\mathcal{ H}}^{(N,K,a)}({\bf D}^2\times {\bf S}^1)\ 
\langle H_1(x,\Omega)\rangle=
\langle U_0(x)\rangle\langle U_1(\Omega)\rangle\\ 
\langle U_1(\Omega)\rangle
 \text{ is invertible}
\end{cases}$$
Here,
for $\epsilon\in\{-1,0,1\}$ we denote by $U_\epsilon$ the unknot with
framing $\epsilon$, and by $H_\epsilon$ the Hopf link 
 with linking number one and both components having
       framing 
$\epsilon$.

A solution of the above is essentially unique.
 In \cite{varso}, we discussed this condition  in
the context of a formal skein theory.

 A framed link $L$ determines by surgery
       a $3$-manifold  denoted by $S^3(L)$ 
(every compact oriented  $3$-manifold can be obtained in this way).
As a consequence of Kirby's theorem \cite{Ki}, if a solution
$\Omega$ exists, then an appropriate normalization
of $\langle L(\Omega,\dots,\Omega)\rangle$
is an invariant of the surgered manifold $M={\bf S}^3(L)$.
It is convenient to choose a solution
$\omega$ of $(K)$ such that
${\langle U_1(\omega)\rangle\langle U_{-1}(\omega)\rangle}=1$.
(It may be necessary to extend $k$ by adding a square root.)
An appropriate normalization is then
$$\tau(M)=\langle U_{1}(\omega)\rangle^{-\sigma(L)}
\langle L(\omega,\dots,\omega)\rangle$$
Here $\sigma(L)$ is the signature of the linking matrix,
i.e. $\sigma(L)=b_+-b_-$, where $b_+$
(resp. $b_-$) is the number of positive (resp. negative)
eigenvalues of the linking matrix $B_L$ associated with $L$.


Recall that we have denoted  by $\mathcal{\Gamma}_{N,K}$
the set of Young diagrams
with at most $N-1$ lines and $K$ columns
(the empty diagram is included).
We set $$\Omega=\sum_{\lambda\in\Gamma_{N,K}}
\langle {\lambda}\rangle\hat y_\lambda$$
Following Yokota \cite{Yo}, we can show the sliding property.
 From this sliding property and
the proposition \ref{reverse} we can deduce
the first part of the Kirby condition
(the handle slide). 
\begin{pro}[Sliding property]\label{pslide}
The Homflypt polynomial of a link in $\mathbf{S}^3$,
which has one of its components  cabled with
 the skein element $\Omega$,
satisfies the equality
in figure \ref{slide}.
\end{pro}
In this figure, the dashed curve means that the component cabled with
$\Omega$ may be non-trivially embedded in the sphere;
the bracket notation for the Homflypt polynomial is omitted.
\begin{figure}
\centerline{
        \begin{pspicture}(0,0)(5,4)
        \psline{->}(3,0)(3,3.5)
        \pscurve[linestyle=dashed]{-}(.6,1.4)(.8,1)(1.5,0.5)(2.2,1)(2.4,1.4)
        \pscurve{->}(2.4,1.4)(2.45,1.7)(2.4,2)
        \pscurve(2.4,2)(2.2,2.5)(1.5,3)(.8,2.5)(.6,2)(.55,1.7)(.6,1.4)
        \put(.8,2.1){$\Omega$}
        \put(4,1.55){\Large$ =$}
        \end{pspicture}
\begin{pspicture}(0,0)(3.5,4)
        \pscurve{->}(.5,1.4)(.45,1.7)(.5,2)(.7,2.6)(1.5,3.1)(2.2,2.8)
(2.8,3)(3,3.5)
        \pscurve(3,0)(2.9,.5)(2.65,.9)(2.5,.9)(2.3,.8)
        \pscurve[linestyle=dashed](.5,1.4)(.7,0.9)(1.5,0.4)(2.3,.8)
        \pscurve[linestyle=dashed]{-}(.6,1.4)(.8,1)(1.5,0.5)(2.2,1)(2.4,1.4)
        \pscurve{-}(2.4,1.4)(2.45,1.7)(2.4,2)
        \pscurve{<-}(2.4,2)(2.2,2.5)(1.5,3)(.8,2.5)(.6,2)(.55,1.7)(.6,1.4)
        \put(.8,2.1){$\Omega$}
        \end{pspicture} 
        }
\caption{\label{slide} Sliding property}
\end{figure}
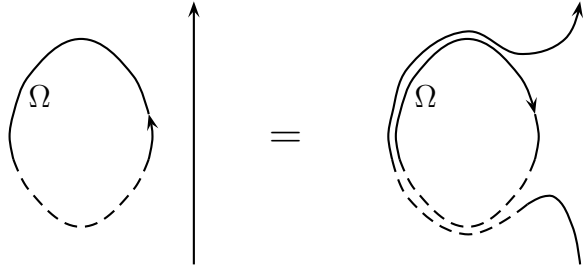
\begin{proof}
Using the definition of $\Omega$, the ramification formula,
and the lemma \ref{rect}.
$$lhs=\sum_\lambda \langle {\lambda } \rangle\sum_{\genfrac
        {}{}{0pt}{2}{\lambda\subset \mu}{|\mu|=|\lambda|+1}}
        \langle \mu\rangle \langle A(\lambda,\mu)\rangle,$$
where $A(\lambda,\mu)$ is the skein element represented in figure
        \ref{proof}.
Note that in the above, we can forget any $\mu$ with $\mu_1=K+1$, because in
        this case the corresponding skein element contains a negligible
morphism.
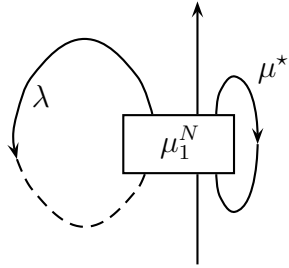
\begin{figure}
\centerline{
        \begin{pspicture}(0,0)(5,4)
        \psline{->}(3,2)(3,3.5)\psline(3,0)(3,1.2)
        \pscurve[linestyle=dashed]{-}(.6,1.4)(.8,1)(1.5,0.5)(2.2,1)(2.3,1.2)
        \pscurve{->}(2.4,2)(2.2,2.5)(1.5,3)(.8,2.5)(.6,2)(.55,1.7)(.6,1.4)
        \put(.8,2.1){$\lambda$}
        \psframe(2,1.2)(3.5,2)
        \pscurve{->}(3.25,2)(3.3,2.3)(3.5,2.5)(3.8,1.6)
        \pscurve(3.8,1.6)(3.5,.7)(3.3,.9)(3.25,1.2)
        \put(2.5,1.5){$\mu_1^N$}
        \put(3.8,2.5){$\mu^\star $}
\end{pspicture}
}\caption{\label{proof} The rectangle trick}
\end{figure}
 We can cut the band colored with 
$\mu_1^N$, move one of the ends along the closed component
 and glue it back again, so that the band
colored with $\mu^\star$ goes along this closed component with the reverse
orientation. Then we
apply  lemma \ref{rect} and the ramification
formula again.
 Using the involution
$\star$ on the set $\Gamma_{N,K}$,
 and also  lemma \ref{qdim},
we obtain the required equality.
\end{proof}

\subsection{The invariant $\tau^{SU({N,K})}$, and
the modular category $\mathrm{H}^{SU(N,K)}$}
Until the end of section \ref{sectionSUN},
we will consider the
choice of the framing parameter $a=s^{\frac{-1}{N}}$.
 More precisely $a$ is a primitive
 $2N(K+N)$-th root of unity, and $s=a^{-N}$.\\
For this choice of parameter $a$, we denote the category
$\overline{\mathrm{H}}^{(N,K,a)}$ by  $\mathrm{H}^{SU(N,K)}$,
and the corresponding  skein functor
by $\mathcal{H}^{SU(N,K)}$.
 \begin{lem}\label{omega}
$$ {\langle U_1(\Omega)\rangle\langle U_{-1}(\Omega)\rangle}
=\langle \Omega \rangle=
(-1)^{\frac{N(N-1)}{2}}
\frac{N(N+K)^{N-1}}{\prod_{1\leq j\leq N}(s^j-s^{-j})^{2(N-j)}}$$
\end{lem}
\begin{proof}
The first equality comes from the sliding property
and  lemma \ref{opc} below.
 For the second equality, we follow Erlijman's
 computation in \cite{Erlij}.
$$\langle \lambda\rangle=\frac{\prod [N+cn]}{\prod [hl]}=
s^{-(N-1)|\lambda|}\mathcal{S}_\lambda(1,s^2,\dots,s^{2(N-1)})$$
Here $\mathcal{S}_\lambda$ is the $\lambda$-indexed Schur symmetric
polynomial
(see ch1 in \cite{Macdo}),
$$\langle \lambda\rangle^2=
\frac{a_{\rho+\lambda}\overline a_{\rho+\lambda}}
{a_{\rho}\overline a_{\rho}}$$
with $\rho=(N-1,N-2,\dots,0)$ and, for
$l=(l_1,\dots,l_N)$
$$a_l=det\left(s^{2(i-1)l_j}\right)
_{1\leq i,j\leq N}$$
We have
$$a_{\rho}\overline a_{\rho}=\prod_{i<j}(s^{2j}-s^{2i})(s^{-2j}-s^{-2i})
=(-1)^{\frac{N(N-1)}{2}}\prod_{\nu=1}^{N-1}(s^\nu-s^{-\nu})^{2(N-\nu)}$$
$$\sum_{\lambda\in \Gamma_{N,K}}
a_{\rho+\lambda}\overline a_{\rho+\lambda}=
\sum_{N+K>l_1>\dots>l_N=0} det\left(s^{2(i-1)l_j}\right)
det\left(s^{-2(i-1)l_j}\right)$$
$$\sum_{\lambda\in \Gamma_{N,K}}
a_{\rho+\lambda}\overline a_{\rho+\lambda}=\frac{1}{(N-1)!}
\sum_{0\leq l_1,\dots,l_{N-1}<N+K}
\sum_{\pi,\pi'\in \mathcal{S}_N} \epsilon_\pi\epsilon_{\pi'}
\prod_i s^{2l_i(\pi(i)-\pi'(i))}$$ 
Here, we have symmetrized the index set for $l$
(and added terms which are zero), and we expanded
the determinants.
$$\sum_{\lambda\in \Gamma_{N,K}}
a_{\rho+\lambda}\overline a_{\rho+\lambda}=\frac{1}{(N-1)!}
\sum_{\pi,\pi'}\epsilon_\pi\epsilon_{\pi'}
\prod_{i=1}^{N-1}\sum_{l_i=0}^{N+K-1}
s^{2l_i(\pi(i)-\pi'(i))}$$
For $\pi\neq\pi'$, a zero term appears,
and the $N!$ remaining terms are all equal. We get
$$\sum_{\lambda\in \Gamma_{N,K}}a_{\rho+\lambda}\overline a_{\rho+\lambda}
=N(N+K)^{N-1}$$
Our formula follows.
 \end{proof}
\begin{lem}\label{opc}
If $\lambda\neq \emptyset$, then the following
morphism is zero in $\mathrm{H}^{SU(N,K)}$.

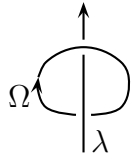
\centerline{
        \begin{pspicture}(0,-0.5)(4,2.5)
        \psline{-}(1.5,0)(1.5,1.3)
        \psline{->}(1.5,1.5)(1.5,2)
        \pscurve{->}(1.4,.5)(1,.6)(.9,.8)(.9,1)
        \pscurve{-}(.9,1)(1,1.2)(1.5,1.4)
        (2,1.2)(2.1,1)(2,.6)(1.6,.5)
        \put(1.6,0){$\lambda$}\put(.5,.6){$\Omega$}
        \end{pspicture}
        }
\end{lem}

\begin{proof}
By the sliding property, for any Young diagram $\mu$ such that
$\lambda\subset\mu$, and $|\mu/\lambda|=1$, we have

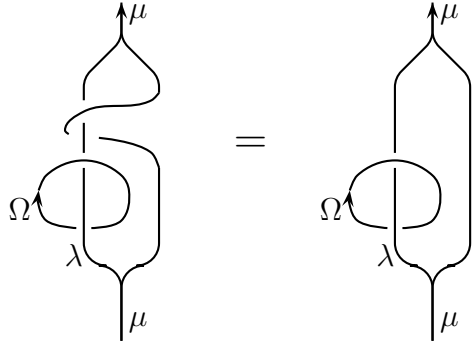
\centerline{
        \begin{pspicture}(0,-1.5)(4,4)
        \psline[linearc=.3]{-}(2,-1)(2,-0,8)(1.5,0)(1.5,1.3)
        \psline[linearc=.3]{<-}(2,3.5)(2,3)(1.5,2.5)(1.5,2.2)
        \psline(1.5,1.5)(1.5,1.9)
        \pscurve{->}(1.4,.5)(1,.6)(.9,.8)(.9,1)
        \pscurve{-}(.9,1)(1,1.2)(1.5,1.4)(2,1.2)(2.1,1)(2,.6)(1.6,.5)
        \put(1.25,0){$\lambda$}\put(.5,.6){$\Omega$}
        \psline[linearc=.3]{-}(2,-1)(2,-0,8)(2.5,0)(2.5,1)(2.5,1.3)
        \psline[linearc=.3](2.5,2.3)(2.5,2.5)(2,3)(2,3.5)
        \pscurve(2.5,1.3)(2.4,1.5)(1.7,1.7)
        \pscurve(1.3,1.75)(1.25,1.8)(1.5,2.05)(2.4,2.2)(2.5,2.3)
        \put(2.1,-.8){$\mu$}\put(2.1,3.3){$\mu$}
\put(3.5,1.5){\Large $=$}
\end{pspicture}
        \begin{pspicture}(0,-1.5)(4,4)
        \psline[linearc=.3]{-}(2,-1)(2,-0,8)(1.5,0)(1.5,1.3)
        \psline[linearc=.3]{<-}(2,3.5)(2,3)(1.5,2.5)(1.5,1.5)
        \pscurve{->}(1.4,.5)(1,.6)(.9,.8)(.9,1)
        \pscurve{-}(.9,1)(1,1.2)(1.5,1.4)(2,1.2)(2.1,1)(2,.6)(1.6,.5)
        \put(1.25,0){$\lambda$}\put(.5,.6){$\Omega$}
        \psline[linearc=.3]{-}(2,-1)(2,-0,8)(2.5,0)(2.5,1)(2.5,2.5)(2,3)(2,3.5)
        \put(2.1,-.8){$\mu$}\put(2.1,3.3){$\mu$}
\end{pspicture}
}
By  proposition \ref{braiding}, if the color $\lambda$ is not killed
by inserting a $0$-framed meridian cabled with  $\Omega$, then
$a^{2|\lambda|}s^{2cn(\mu/\lambda)}=1$  for any
$\mu$ as above with $\langle\mu\rangle \neq 0$.
If two such $\mu$ exist, with added cells of
respective indices $(i,j)$ and $(i',j')$, then
$s^{2(j'-j+i-i')}=1$ and $-(N+K)<j'-j+i-i'<N+K$. This implies
 $(i,j)=(i',j')$.
The remaining diagrams are those with $K$ columns and
$n$ lines, $n<N$.
In this case we get $a^{2nK}s^{-2n}=a^{2n(K+N)}=1$.
 The order of $a$ implies $n=0$.
\end{proof}

We suppose that $N(N+K)$ is invertible. We set
$\eta^{-2}=\langle \Omega \rangle$ (we extend $k$ if necessary),
$\omega=\eta\Omega$, $\Delta=\langle U_{1}(\omega)\rangle$.
\begin{thm}
There exists an invariant of compact oriented $3$-manifolds 
defined on a surgery presentation by the following formula. 
$$\tau^{SU({N,K})}(\mathbf{S}^3(L))=\Delta^{-\sigma(L)}
\langle L(\omega,\dots,\omega)\rangle
$$
\end{thm}
This invariant can be extended to manifolds
with (colored) links by the formula
$$\tau^{SU(N,K)}(\mathbf{S}^3(L),K)=\Delta^{-\sigma(L)}
\langle L(\omega,\dots,\omega)\cup K\rangle
$$
The above is a consequence of the Kirby theorem.
 Using \cite{Turaev}, the following  constructs the associated TQFT
(and proves again the theorem above).
\begin{thm}
The category $\mathrm{H}^{SU(N,K)}$ is a modular category
with $ \Gamma_{N,K}$ as a representative set of isomorphism classes of
 simple objects.
\end{thm}
\begin{proof}
We have to check invertibility of the $S$ matrix whose
entries, indexed by $\Gamma_{N,K}$, are
the evaluations of the Homflypt polynomial
for a Hopf link $H_0$ ($0$-framed, and with linking $+1$),
whose components are colored with the corresponding indices.
The lemma below gives the result (the matrix $\overline S$
is the conjugate of $S$).
\end{proof}
\begin{lem}
One has ($I$ is the unit matrix)
$$S\overline S=
\langle \Omega\rangle I$$
\end{lem}
\begin{proof}
The $(\lambda,\mu)$-indexed entry of the matrix
$S\overline S$ can be written
$$u_{\lambda\mu}=\langle H_0(\hat y_\lambda  \hat y_{\mu^*}),\Omega)\rangle$$
 We write $u_{\lambda\mu}$ as the quantum trace of the morphism
 $\gamma_{\lambda\mu}$ represented by a
$\lambda\otimes\mu^*$-colored band together
with an $\Omega$-cabled meridian.
 Using that the objects in $\Gamma_{N,K}$
dominate, and the lemma \ref{opc}, we obtain
that there exists a finite family $(\alpha_i,\beta_i)$,
$\alpha_i\in  {\mathrm{H}}^{SU(N,K)}(\lambda\otimes\mu^*,0)$,
$\beta_i\in {\mathrm{H}}^{SU(N,K)}(0,\lambda\otimes\mu^*)$,
such that
$\gamma_{\lambda\mu}=\langle\Omega\rangle\sum_i{ \alpha_i \beta_i}$.\\
If $\lambda$ and $\mu$ are distinct in $\Gamma_{N,K}$, then
 the two modules above are zero, hence we have $u_{\lambda\mu}=0$.\\
 If $\lambda=\mu$, then this two modules are generated by the duality
 morphisms. This gives
$\sum_i\alpha_i \beta_i=\zeta d_{\lambda^*}b_\lambda$.
 We obtain $\zeta=\frac{1}{\langle \lambda \rangle}$ from
 the product $b_\lambda\gamma_{\lambda\lambda}d_{\lambda^*}=
\langle \lambda \rangle\langle \Omega\rangle$.
 We can conclude that $u_{\lambda\lambda}=\langle \Omega\rangle$.
\end{proof}
As already stated, using Turaev's work, modularity
gives the TQFT.  The universal construction
of \cite{BHMV3} could also be applied here.
The normalized invariant of connected closed $3$-manifolds $M$ equipped
with $p_1$-structure (or $2$-framing) $\alpha$ is defined
by $$Z(M,\alpha)=\eta \kappa^{-\sigma(\alpha)}
\tau^{SU(N,K)}(M)$$
with 
$\kappa^3=\Delta$.
We obtain a TQFT functor $(V,Z)$ on the cobordism category
$C^2_{p_1}$ of $p_1$-surfaces and (equivalence classes
of) cobordisms.

Unfortunately the description of structured surfaces
in \cite{Turaev} and \cite{BHMV3} are not the same, however
this has essentially no influence on the description
of the TQFT modules.
 
The vectors $Z(\mathbf{D}^2\times \mathbf{S}^1,\hat y_\lambda)$,
$\lambda\in \Gamma_{N,K}$ form a basis for the TQFT
module $V(\mathbf{S}^1\times \mathbf{S}^1)$.
Moreover this basis is orthonormal with respect
to the natural hermitian form on this module.

\noindent{\bf The fusion algebra and Verlinde dimension formula.}
The algebra structure on the skein module of the solid torus
induces the fusion algebra structure on
$V(\mathbf{S}^1\times \mathbf{S}^1)$.
 From modularity, we get the structure constants
 $$\hat y_\lambda \hat y_\mu=\sum_\nu c_{\lambda\mu}^\nu \hat y_\nu$$
 with $c_{\lambda\mu}^\nu$ equal to the rank of the
 module associated to a sphere, with two incoming points and
 one outgoing point
 colored respectively by $\lambda$, $\mu$ and $\nu$.
 A combinatorial description for these ranks should be obtainable
 using  \cite{Good}.

 For a genus $g$ closed surface $\Sigma_g$, we can compute the dimension
 $d_g$ of the TQFT module $V(\Sigma_g)$, which is equal
 to the invariant $Z(\Sigma_g\times \mathbf{S}^1)$.
  The result is given by the Verlinde formula.
  $$d_g=
  \left((N+K)^{(N-1)}N\right)^{g-1}
  \sum_{N+K>l_1>l_2>\dots>l_N=0}\ \prod_{1\leq i<j\leq N}
  \left(\frac{-1}{(s^{l_j-l_i}-s^{l_i-l_j})^2}\right)^{g-1}$$
  
\subsection{The invariant $\tau^{PSU(N,K)}$ and the
modular category $\mathrm{H}^{PSU(N,K)}$}
The algebra ${\mathcal{H}}^{SU(N,K)}(\mathbf{D}^2\times\mathbf{S}^1)$
is $N$-graded.\\
Set $\Gamma_{N,K}^0=\{\lambda\in  \Gamma_{N,K},\ |\lambda|
  \equiv 0\ \mathrm{mod}\ N\}$.
  
We can see that
  $\Omega_0=\sum_{\lambda\in\Gamma_{N,K}^0}
\langle \lambda\rangle\hat y_\lambda$
  satisfies the handle slide condition (an arc with a $0$-graded color
  can slide over a component cabled with $\Omega_0$).
\begin{lem}
If $d=\mathrm{gcd}(N,K)$ is even, and $N'=\frac{N}{d}$ $K'=\frac{K}{d}$
are both odd, then
$$\langle U_1(\Omega_0)\rangle=0\ .$$
In all other cases one has that
$$ {\langle U_1(\Omega_0)\rangle\langle U_{-1}(\Omega_0)\rangle}
=\langle \Omega_0 \rangle=
(-1)^{\frac{N(N-1)}{2}}
\frac{(N+K)^{N-1}}{\prod_{ j=1}^{N-1}(s^j-s^{-j})^{2(N-j)}}$$
\end{lem}
We say that the rank-level $(N,K)$ is {\it spin} if
$d=\mathrm{gcd}(N,K)$ is even, and $N'=\frac{N}{d}$, $K'=\frac{K}{d}$
are both odd. The terminology will be justified
in section \ref{refined}.
 
We suppose now that the rank-level $(N,K)$ {\em is not spin}, and 
that $N+K$ is invertible.

We set
$\eta_0^{-2}=\langle \Omega_0 \rangle$ (we extend $k$ if necessary),
$\omega_0=\eta_0\Omega_0$,
$\Delta=\langle U_{1}(\omega)\rangle$
(we can choose the sign of $\eta_0$
so that this $\Delta$ is the same as above).
\begin{thm}
There exists an invariant of compact oriented $3$-manifolds
defined on a surgery presentation by the following formula. 
$$\tau^{PSU(N,K)}(\mathbf{S}^3(L))=\Delta^{-\sigma(L)}
\langle L(\omega_0,\dots,\omega_0)\rangle
$$
\end{thm}
\begin{rem}
 For $\mathrm{gcd}(N,K)=1$ this invariant and an underlying
modular category obtained from the quantum group are known
(see \cite{MW}).
\end{rem}
\begin{rem}
There exists a refined invariant $\tau^{SU(N,K)}(M,c)$,
with $c$ a cohomology class in $H^1(M,\mathbb{Z}/gcd(N,K))$ such that
$\tau^{PSU(N,K)}(M)=\tau^{SU(N,K)}(M,0)$ (see section
\ref{refined}).\end{rem}
\noindent
{\bf The modular category.}
The objects of the category $\mathrm{H}^{SU(N,K)}$
are graded by the algebraic 
number of points in their expansion (signed with the orientation).
Let $\mathrm{H}^{PSU(N,K)}$ be the full
subcategory of $\mathrm{H}^{SU(N,K)}$, whose
objects are zero graded modulo $N$.
 We can show the following (we give no details).
 \begin{thm}
The category $\mathrm{H}^{PSU(N,K)}$ is a modular category
with $\Gamma_{N,K}^0$ as a representative set of isomorphism classes of
 simple objects.
 \end{thm}

\section{The modular category $\widetilde{\mathrm H}^{N,K}$
and the  invariant $\tilde\tau_{N,K}$}
\label{reduced}
In this section, we work at rank $N$ as before
(i.e. $v=s^{-N}$), but level
$K$ will mean that, in the integral domain $k$,
$$\begin{cases}
s \text{ has order $2(N+K)$ if $N+K$ is even,}\\
s \text{ has order $N+K$ if $N+K$ is odd.}
\end{cases}$$
In both cases $s^2$ has order $N+K$; we have $v=\varepsilon s^K$,
with $\varepsilon=(-1)^{N+K+1}$.\\
Let $d=\mathrm{gcd}(N,K)$, $N=dN'$, $K=dK'$.
Motivated by the formulas in proposition
\ref{braiding},
we would like to fix the framing parameter $a$ in such a way that
the order of the multiplicative subgroup
generated by $a^Ns$ and $a^{K}s^{-1}$ is as small as possible.
 Note that this order is at least $d$ in the odd case, and $2d$
in the even case. We show that this lower bound can be realized.\\
To simplify the discussion
we suppose in the case $N+K$ even that 
$N'$ is odd;
we set $d=\alpha \beta$ with
$\mathrm{gcd}(\alpha,2K')=\mathrm{gcd}(\beta,N')=1$.
 Recall that
 $\varepsilon=(-1)^{N+K+1}$.
\begin{lem}\label{res}
We can choose
the framing parameter $a$ so that $(a^Ns)^\alpha=1$ and
$(a^Ks^{-1})^\beta=\varepsilon$.
\end{lem}
\noindent (It may be necessary to extend the scalars.)
\begin{proof}
The problem is to find a common solution to certain polynomial equations.
Here the resultant is
$$\text{Res}_{a}(a^{\alpha N}-s^{-\alpha},
a^{\beta K}-\varepsilon s^{\beta})=
(s^{-K}-\varepsilon s^ {N})^d$$
\end{proof}
\begin{rem}
If $N+K$ and $N'$ are both even, then we set
$d=\alpha\beta$, with
$\mathrm{gcd}(\alpha,K')=\mathrm{gcd}(\beta,N')=1$,
and in the above lemma we require that
$(a^Ns)^\alpha=-1$ and
$(a^Ks^{-1})^\beta=1$.
 Following \cite{KT2}, we can show a level-rank duality
 formula, and recover this case by exchanging
 $N$ and $K$.
 \end{rem}
\noindent{\bf In this section, and in section \ref{refined},
we suppose that the framing parameter $a$
satisfies the condition of lemma \ref{res}.}
Recall that the category $\mathrm{H}^{(N,K,a)}$
is the quotient of the $\mathcal{C}^{N,K}$-extended
Hecke category by negligible morphisms.
 We define the  category $\widetilde{\mathrm H}^{N,K}$
 as follows. Objects are those of $\mathrm{H}^{(N,K,a)}$,
 and the modules of morphisms are obtained from
 those of $\mathrm{H}^{(N,K,a)}$ by adding as generators
 colored ribbon graphs in which incoming or outgoing
$1$-valent vertices colored with $(1^N)^{\otimes\alpha}$
or $(K)^{\otimes\beta}$ are allowed,
and quotienting by the relations given by gluing
an incoming $(1^N)^{\otimes\alpha}$-colored 
(resp. $(K)^{\otimes\beta}$-colored) vertex
with an outgoing one.
 The corresponding skein functor is denoted by
 $\widetilde{\mathcal{H}}^{N,K}$.\\[5pt]
\noindent{\em Exercise.} a) Let $\mu\in \overline{\Gamma}_{N,K}$ be a
Young diagram of the form $\mu=1^N+\lambda$. Show that the objects $\mu$
and $1^N\otimes \lambda$ are isomorphic in the category
$\widetilde{\mathrm{H}}^{N,K}$.\\
b) Show that the objects $(1^N)^{\otimes K}$
and $(K)^{\otimes N}$ are isomorphic in the category
$\widetilde{\mathrm{H}}^{N,K}$.\\
c) Show that the Homflypt invariant extends to an isomorphism
$\widetilde{\mathcal{H}}^{N,K}({\bf S}^3)\simeq k$.

\vspace{5pt}
Set $\dot\Gamma_{N,K}=\{(1^N)^{\otimes i}\otimes \lambda,\ 0\leq i<\alpha
\text{ and }
\lambda \in \Gamma_{N,K}\}$.\\
The category $\widetilde{\mathrm{H}}^{N,K}$ is a  ribbon category,
and $\dot\Gamma_{N,K}$ is a  finite set of dominating simple objects.
  We have that  $\widetilde{\mathrm H}^{N,K}$ is a pre-modular category.
This uses the involution $\star$ on the set $\dot\Gamma_{N,K}$
defined as follows.

For $x=(1^N)^{\otimes i}\otimes \lambda\in\dot\Gamma_{N,K}$,
we set $x^\star=\lambda^\star\otimes (1^N)^{\otimes i'}$,
where $i'\in\{0,\dots,\alpha-1\}$ is such
that the number of points in $x\otimes x^\star$
is a multiple of $N\alpha$.

We can now proceed similarly as we did in section
\ref{sectionSUN}.
\begin{pro}
Let $L$ be a framed link in the $3$-sphere, and let $L'$ be
the link obtained from
 $L$ by reversing the orientation of the first component.
 Then
for any $x\in \dot{\Gamma}_{N,K}$, one has
$$\langle L(\hat x,\dots)\rangle=
\langle L'(\hat x^\star,\dots)\rangle$$
\end{pro}
Set $$\dot \Omega=\sum_{x\in\dot\Gamma_{N,K}}\langle x\rangle\hat1_x=
\sum_{i=0}^{\alpha-1}
\hat g_N^{i}\Omega$$
where $\Omega$ is defined as before by
$$\Omega=\sum_{\lambda\in \Gamma_{N,K}} \langle \lambda\rangle
\hat y_\lambda$$
\begin{pro}[Sliding property]
The Homflypt polynomial of a link in $\mathbf{S}^3$,
which has one of its components  cabled with
the skein element $\dot\Omega$,
satisfies the equality
in figure \ref{slide} (with $\Omega$ replaced by $\dot\Omega$).
\end{pro}
There is an action of $\mathbb{Z}/N$ on the set
$\dot \Gamma_{N,K}$ defined as follows.
The generator of $\mathbb{Z}/N$ acts by
$$
(1^N) ^{\otimes i}\otimes\lambda\mapsto
(1^N) ^{\otimes( i'+\lambda_{N-1})}\otimes
((K,\lambda)-\lambda_{N-1}^N)\ ,$$
where $i'\in\{0,\dots,\alpha-1\}$ is congruent to
$i+\lambda_{N-1}$ modulo $\alpha$.\\
The idea is that we add to the diagram a line with
$K$ cells, and then each column which has  $N$ cells
is replaced by an added copy of $1^N$.

Restricting this action to the group generated by $\beta$,
we obtain an action of $\mathbb{Z}/{\alpha N'}$.
 By considering the degree mod $\alpha N$, we can see
 that this action is free.
We denote by $\widetilde \Gamma_{N,K}$ a subset of
$\dot\Gamma_{N,K}$ which is a representative set for the orbits.
 We set $$\tilde \Omega=\sum_{u\in \widetilde \Gamma_{N,K}}
\langle u\rangle {\hat1_u}.$$
Cabling $\dot \Omega$ gives the same result as
cabling $N'\alpha \tilde \Omega$,
so that $\tilde \Omega$ also satisfies the
 sliding property.
 \begin{lem}\label{tildeomega}
$$ {\langle U_1(\tilde\Omega)\rangle\langle U_{-1}(\tilde\Omega)\rangle}
=\langle \tilde\Omega \rangle=
(-1)^{\frac{N(N-1)}{2}}
\frac{d(N+K)^{N-1}}{\prod_{j=1}^{N-1}(s^j-s^{-j})^{2(N-j)}}$$
\end{lem}
We suppose that $(N+K)$ is invertible. We set
$\tilde{\eta}^{-2}=\langle \tilde \Omega \rangle$
(we extend $k$ if necessary),
$\tilde\omega=\tilde\eta\tilde\Omega$,
$\delta=\langle U_{1}(\tilde\omega)\rangle$.
\begin{thm}
a) There exists an invariant of compact oriented $3$-manifolds 
defined on a surgery presentation by the following formula. 
$$\tilde\tau_{N,K}(\mathbf{S}^3(L))=\delta^{-\sigma(L)}
\langle L(\tilde\omega,\dots,\tilde\omega)\rangle
$$
b) The category $\widetilde{\mathrm H}^{N,K}$ is a modular category
with $\widetilde \Gamma_{N,K}$ as a representative set
of isomorphism classes of
 simple objects.
\end{thm}

Now we give the relation between the invariant
$\tau^{SU({N,K})}$ and $\tilde \tau_{N,K}$ in the following
reduction theorem.
Here $\tau^{U(1)}(M,\zeta)$ is a version of
the invariant derived from linking
matrices in \cite{MOO} ($U(1)$ invariant) for a root
of unity $\zeta$, whose
order is $N'$ (resp. $2N'$) if $N'$ is odd
(resp. if $N'$ is even).
\begin{thm}[Reduction formula]\label{reduction}
For every manifold $M$,
one has
$$\tau^{SU({N,K})}(M)=\tau^{U(1)}(M,\zeta)\tilde\tau_{N,K}(M)$$
\end{thm}
We will give the proof in  section \ref{last}.
 We write this equality in the ring $k$
 where the  invariant $\tau^{SU(N,K)}$ has been  defined.
 The definition of the reduced invariant
$\tilde\tau_{N,K}$ requires a choice of the parameters
in $k$; we denote by $(\tilde s,\tilde v, \tilde a)$
this choice.
\begin{verse}
 If $d$ is even, then
$\tilde s=s$, $\tilde v=v$ and
$\zeta=(a^Ks^{-1})^{K\beta^2}$
(the case $N'$ even is not excluded).\\
If $d$ is odd, then $\tilde s=-s$, $\tilde v=-v$
and $\zeta=((-a)^Ks^{-1})^{K\beta^2}$.
\end{verse}
 The $U(1)$ invariant in the formula above is then defined by
$$\tau^{U(1)}(\mathbf{S}^3(L),\zeta)=
\left(\frac{\Delta}{\delta}\right)^{-\sigma(L)}
\left(\frac{\eta}{\tilde \eta}\right)^m
\sum_{j\in(\mathbb{Z}/{N'})^m}
\zeta^{{}^tjB_L j} .$$
 
\section{Refined invariants}\label{refined}
\subsection{Spin structures modulo an even integer}
 In \cite{Mu}, H. Murakami
 stated a decomposition formula for the invariant
$\tau^{SU({N,K})}$ using  some spin type structures (see remark 2.7
in his paper).
He observed that for $N=2$ these are spin structures,
and the corresponding refinements were studied in \cite{KM} and \cite{Bl1}.
For $N>2$ he only gave a combinatorial description of the structures,
and asked for a topological interpretation.
We recall here the topological definition for these structures
which we gave in \cite{varso}.

Suppose $d$ is an even integer. Then there exists, up to homotopy, a
unique non trivial map
$g: BSO\rightarrow K(\mathbb{Z}/d,2)$.
Define the fibration
$$\pi_d: BSpin(\mathbb{Z}/d)\rightarrow BSO$$
to be the pull-back, using $g$, of the path
fibration over $K(\mathbb{Z}/d,2)$.
 The space $BSpin(\mathbb{Z}/d)$ is a classifying space
for the non trivial central extension
of the Lie group $SO$ by $\mathbb{Z}/d$,
which we denote by $Spin(\mathbb{Z}/d)$.

Now we can use the fibration $\pi_d$ to define structures
(see [St]).  Let $\gamma_{Spin(\mathbb{Z}/d)}=\pi_d^*(\gamma_{SO})$
be the pull-back of the canonical vector bundle over $BSO$.

\begin{dfn} A   spin structure 
 with mod $d$ coefficients on a manifold $M$
is an homotopy class
of fiber maps from the stable tangent bundle $\tau_M$ to
$\gamma_{Spin(\mathbb{Z}/d)}$.
\end{dfn}
If non empty, the set of these structures, denoted
by $Spin(M;\mathbb{Z}/d)$, is affinely isomorphic to
$H^1(M;\mathbb{Z}/d)$, by obstruction theory.
 Moreover the obstruction for existence is a class
 $w_2(M;\mathbb{Z}/d)\in H^2(M;\mathbb{Z}/d)$, which is the image of
 the Stiefel-Whitney class $w_2(M)$ by the homomorphism
 induced by the inclusion of coefficients
 $\mathbb{Z}/2\hookrightarrow \mathbb{Z}/d$.
 As the Stiefel-Whitney class $w_2(M)$ is zero for every compact oriented
$3$-manifold, this shows that
 spin  structures modulo $d$ exist on every
$3$-manifold $M$.
The following theorem gives a combinatorial description
for these structures.  Recall that a surgered manifold
$M=\mathbf{S}^3(L)$ is the boundary of a
$4$-manifold $W_L$ called the trace of the surgery.
To each $\sigma\in Spin(M;{\mathbb{Z}/d})$ is associated
a relative obstruction $w_2(\sigma;{\mathbb{Z}/d})$
in $H^2(W_L,M;\mathbb{Z}/d)$.
 The group $H^2(W_L,M;\mathbb{Z}/d)$ is free of rank $m=\sharp L$.
 Taking the coordinates of the relative obstruction
 we get a map $\psi_L: Spin(M;{\mathbb{Z}/d})
\rightarrow\left(\mathbb{Z}/d\right)^m$.

\begin{thm} The map $\psi_L$ is injective, and its image
is the set of those $(c_1,\dots,c_m)$ which are
solutions of the following $(\mathbb{Z}/d)$-characteristic
equation
$$B_L\left(\begin{array}{c} c_1\\ \vdots\\ c_m\end{array}\right)=
{\frac{d}{2}}\left(\begin{array}{c} b_{11}\\ \vdots\\
 b_{mm}\end{array}\right)\ \ ({\rm mod}
\ d) .$$
{\rm Here the $b_{ii}$ are the diagonal values of
the linking matrix $B_L$.}
\end{thm}
\begin{proof} First we compute the absolute obstruction
$w_2(W_L;\mathbb{Z}/d)=\xi_*(w_2(W_L))$, where $\xi_*$ is induced by the
morphism
of coefficients $\xi: \mathbb{Z}/2\hookrightarrow \mathbb{Z}/d$.
If $x$ is an integral $2$-cycle in $W_L$ with self-intersection $x.x$ 
and $[x]_\nu$ denotes its homology class
modulo an integer $\nu$, $w_2(W_L)\in H^2(W_L;\mathbb{Z}/2)$
is determined by the equation
$$\forall x\ <w_2(W_L),[x]_2>=x.x\ ({\rm mod}\ 2) .$$
Hence $w_2(W_L;\mathbb{Z}/d)\in H^2(W_L;\mathbb{Z}/d)$
is determined by 
$$\forall x\ <w_2(W_L;\mathbb{Z}/d),[x]_d>=\xi(x.x)=\frac{d}{2}x.x
\ ({\rm mod}\ d) .$$
Now by functoriality, the relative obstruction lives
in the inverse image of the absolute one under
the map induced by inclusion
$H^2(W_L,M;\mathbb{Z}/d)\rightarrow H^2(W_L;\mathbb{Z}/d)$.
Using the affine structure over $H^1(M;\mathbb{Z}/d)$,
we obtain an affine bijection between $Spin(M;\mathbb{Z}/d)$
and this inverse image.  Whence we have the lemma
by writing the equation above using coordinates.
\end{proof}  
 
There is a formula for the bijection
$\psi_{L,L'}$ corresponding to a Kirby move.
 Using the $\mathbb{Z}/d$-characteristic
equation we see that the coefficient for a trivial component
with framing $\pm 1$ is $d/2$.
 For the usual positive
Fenn-Rourke move, the formula is
$$\psi_{L,L'}(c_1,\dots,c_{m-1},d/2)
=(c_1,\dots,c_{m-1},c'_m)$$
$${\rm with}\ \  c'_m=d/2-\sum_i b'_{im}c_i .$$
Here $b'_{im}$ is the $(i,m)$-indexed
coefficient of the matrix $B_{L'}$. 
\subsection{Spin refinements}
We consider here the reduced theory
of section \ref{reduced} in the spin case. Recall
that this means
$$\begin{cases}
d=\mathrm{gcd}(N,K)\text{ is even},\\
N'=\frac{N}{d}\text{ and } K'=\frac{K}{d}\text{ are odd}.
\end{cases}$$

We decompose the skein element
$\tilde\omega=\sum_\nu \tilde\omega_\nu$
 according
to the $\mathbb{Z}/d$-grading of the algebra
$\widetilde{\mathcal{H}}^{N,K}({\mathbf{D}^2\times\mathbf{S}^1})=
\bigoplus_\nu \widetilde{\mathcal{H}}^{N,K}_\nu
({\mathbf{D}^2\times\mathbf{S}^1})$.

\begin{thm}
Provided $c=(c_1,\dots,c_m)$ satisfies the modulo $d$
characteristic condition,
the formula
$$\tilde\tau_{N,K}^{\mathrm{spin}}(M,\sigma)=
\delta^{-\sigma(L)}{\langle L(\tilde\omega_{c_1},
\dots,\tilde\omega_{c_m})\rangle}
$$
defines an invariant of the surgered manifold $M={\bf S}^3(L)$
equipped with the modulo $d$ spin  structure
$\sigma=\psi_L^{-1}(c_1,\dots,c_m)$
.\\
Moreover,
$$\forall M\ \ \tilde\tau_{N,K}(M)=\sum_{\sigma\in Spin(M;\mathbb{Z}/d)}
\tilde\tau_{N,K}^{\mathrm{spin}}(M,\sigma) .$$
\end{thm}
\begin{proof}
The following graded version of the sliding property can be
derived from the proof of proposition \ref{pslide}.
\begin{lem}[Graded sliding property]
 The Homflypt polynomial of a link in $\mathbf{S}^3$,
which has one of its components  cabled with
 the skein element $\tilde\omega_\nu$,
satisfies the equality
in figure \ref{gslide}.
\end{lem}
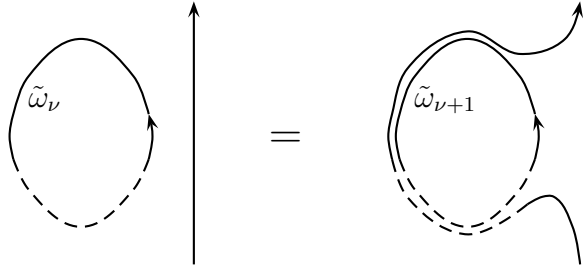
\begin{figure}
\centerline{
        \begin{pspicture}(0,0)(5,4)
        \psline{->}(3,0)(3,3.5)
        \pscurve[linestyle=dashed]{-}(.6,1.4)(.8,1)(1.5,0.5)(2.2,1)(2.4,1.4)
        \pscurve{->}(2.4,1.4)(2.45,1.7)(2.4,2)
        \pscurve(2.4,2)(2.2,2.5)(1.5,3)(.8,2.5)(.6,2)(.55,1.7)(.6,1.4)
        \put(.8,2.1){$\tilde\omega_\nu$}
        \put(4,1.55){\Large$ =$}
        \end{pspicture}
\begin{pspicture}(0,0)(3.5,4)
        \pscurve{->}(.5,1.4)(.45,1.7)(.5,2)(.7,2.6)(1.5,3.1)(2.2,2.8)
(2.8,3)(3,3.5)
        \pscurve(3,0)(2.9,.5)(2.65,.9)(2.5,.9)(2.3,.8)
        \pscurve[linestyle=dashed](.5,1.4)(.7,0.9)(1.5,0.4)(2.3,.8)
        \pscurve[linestyle=dashed]{-}(.6,1.4)(.8,1)(1.5,0.5)(2.2,1)(2.4,1.4)
        \pscurve{->}(2.4,1.4)(2.45,1.7)(2.4,2)
        \pscurve{-}(2.4,2)(2.2,2.5)(1.5,3)(.8,2.5)(.6,2)(.55,1.7)(.6,1.4)
        \put(.8,2.1){$\tilde\omega_{\nu+1}$}
        \end{pspicture} 
        }
        \caption{\label{gslide} Graded sliding property}
\end{figure}
Using this sliding property, we get
$$\forall \nu \ \forall x_\nu\in
\widetilde{\mathcal{H}}^{N,K}_\nu(\mathbf{D}^2\times\mathbf{S}^1)\ 
\langle H_1(x_\nu,\tilde\omega_{d/2-\nu})\rangle=
\langle U_0(x_\nu)\rangle\langle U_1(\tilde\omega_{d/2})\rangle$$
This shows that the spin handle slide condition
is satisfied.

We can deduce from the lemma
\ref{declem} that $\langle U_\epsilon(\tilde\omega_{d/2})\rangle=
\langle U_\epsilon(\tilde\omega)\rangle$ is invertible, and hence
we have  that the invariant $\tilde\tau_{N,K}^{\mathrm{spin}}$
is well defined.

Now we show the decomposition formula.
We can write
$$\langle L(\tilde\omega,
\dots,\tilde\omega)\rangle=
\sum_c\langle L(\tilde\omega_{c_1},
\dots,\tilde\omega_{c_m})\rangle .$$
The result is contained in the lemma \ref{charac}.
\end{proof}

\begin{lem}\label{charac}
If $c$ does not satisfy the modulo $d$
characteristic condition, then
$$\langle L(\tilde\omega_{c_1},
\dots,\tilde\omega_{c_m})\rangle=0 .$$
\end{lem}
\begin{proof}
Let $L=(L_1,\dots,L_m)$. Up to a permutation of the components,
we have to show that
$\langle L(\tilde\omega_{c_1},
\dots,\tilde\omega_{c_m})\rangle=0$,
if $\sum_{j=1}^m b_{1j} c_j\neq \frac{d}{2}\text{ mod }d$.
 The proof is in three steps.\\
 If $L_1$ is a $1$-framed unknot, then, using the sliding property,
we get the result from the  lemma \ref{declem}.\\
If $L_1$ is an unknot with any framing, note that (up to an invertible
factor) we can add to the link a trivial $\pm 1$-framed $U_\epsilon$
cabled with $\omega_{d/2}$; we then use the sliding property
to move this unknot around $L_1$, and this adds $\epsilon$
to the framing.  This allows us to reduce the problem to the
preceding case.\\
In the general case, the component $L_1$ can be unknotted
by changing some crossings, and inserting a $\pm 1$-framed unknot cabled
with $\omega_{d/2}$ around the crossing, in such a way that
its linking number with $L_1$ is zero.
By the sliding property, this amounts to multiplying by an invertible
element, whence we have the result by using the above.
\end{proof}
\begin{lem}\label{declem}
For $\epsilon=\pm1$,
$$\langle U_\epsilon(\tilde\omega_\nu)\rangle=
0\text{ if }\nu\neq \frac{d}{2}.$$
\end{lem}
\begin{proof}
We compute $\langle U_1(\tilde \omega_\nu\rangle
\langle U_{-1}(\tilde \omega_\nu)\rangle$.
 The graded sliding property shows that this is
equal to the Homflypt invariant of the following colored
link.

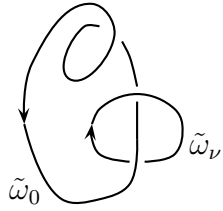
\centerline{
        \begin{pspicture}(-1,-0.5)(4,3)
        \psline{-}(1.5,.8)(1.5,1.3)
        \pscurve{-}(1.5,1.5)(1.4,1.9)(1.3,2.1)
        \pscurve{->}(1,2.3)(.8,2.3)(.6,1.6)(1.2,2.2)
        (.9,2.6)(.1,1.8)(0,1)
        \pscurve{-}(0,1)(.5,0)(1,0)(1.4,.2)(1.5,.8)
        \pscurve{->}(1.4,.5)(1,.6)(.9,.8)(.9,1)
        \pscurve{-}(.9,1)(1,1.2)(1.5,1.4)
        (2,1.2)(2.1,1)(2,.6)(1.6,.5)
        \put(-.2,0){$\tilde\omega_0$}
\put(2.2,.6){$\tilde\omega_\nu$}
        \end{pspicture}
}
        
Using the vanishing lemma \ref{gvanish} below,
 whose proof is adapted from the one given in
\ref{opc}, we obtain
$$\langle U_1(\tilde \omega_\nu\rangle
\langle U_{-1}(\tilde \omega_\nu)\rangle=
\tilde\eta\sum_{i=0}^{\alpha-1}
(a^Ns)^{2\nu i}
\sum_{j=0}^{\beta-1}(-1)^j(a^Ks^{-1})^{2\nu j}.$$
If $\langle U_1(\tilde \omega_\nu\rangle
\langle U_{-1}(\tilde \omega_\nu)\rangle$ is not zero, then
we have
$\nu\equiv 0$ mod $\alpha$, and $\nu\equiv \frac{\beta}{2}$ mod $\beta$.
 \end{proof}
\begin{lem}\label{gvanish}
If $\lambda$ is a Young diagram
in $\Gamma_{N,K}\setminus\{K^j,0\leq j<N\}$, then, for any
$\nu\in\mathbb{Z}/d$, the following 
morphism is zero in $\widetilde{\mathrm{H}}^{N,K}$.

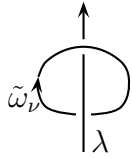
\centerline{
        \begin{pspicture}(0,-0.5)(4,2.5)
        \psline{-}(1.5,0)(1.5,1.3)
        \psline{->}(1.5,1.5)(1.5,2)
        \pscurve{->}(1.4,.5)(1,.6)(.9,.8)(.9,1)
        \pscurve{-}(.9,1)(1,1.2)(1.5,1.4)
        (2,1.2)(2.1,1)(2,.6)(1.6,.5)
        \put(1.6,0){$\lambda$}\put(.5,.6){$\tilde\omega_\nu$}
        \end{pspicture}
        }
\end{lem} 
\subsection{Cohomological refinements}
If the rank-level is not spin, we can proceed
similarly. This time we have that
$\langle U_1(\tilde\omega_\nu)\rangle=0$, unless
$\nu=0\text{ mod }d$.
\begin{thm}
Provided $c=(c_1,\dots,c_m)$ is in the kernel
of the linking matrix modulo $d$,
the formula
$$\tilde\tau_{N,K}^{\mathrm{coho}}(M,\sigma)=
\delta^{-\sigma(L)}{\langle L(\tilde\omega_{c_1},
\dots,\tilde\omega_{c_m})\rangle}
$$
defines an invariant of the surgered manifold $M={\bf S}^3(L)$
equipped with the cohomological
$\sigma\in H^1(M,\mathbb{Z}/d)$ corresponding to $c$.\\
Moreover,
$$\forall M\ \ \tilde\tau_{N,K}(M)=\sum_{\sigma\in H^1(M;\mathbb{Z}/d)}
\tilde\tau_{N,K}^{\mathrm{coho}}(M,\sigma)$$
\end{thm}
\section{Proof of the reduction formula}\label{last}

We show the reduction formula (theorem \ref{reduction}) in the
spin case. 
We proceed as follows.
 We  use the modulo $\beta$ grading to construct, as we did for
 $\tilde \tau_{N,K}$, a spin refinement of
the invariant $\tau^{SU(N,K)}$ satisfying the
decomposition formula
$$\tau^{SU({N,K})}(M)=
\sum_{\sigma\in Spin(M,\mathbb{Z}/\beta)}\tau^{SU({N,K})}(M,\sigma)$$
The same can be done with the reduced invariant
$\tilde\tau_{N,K}$.
 Note that here we consider the modulo $\beta$ grading,
so that the decomposition formula below is not that
of the preceding section if $\alpha\neq1$.
$$\tilde\tau_{N,K}(M)=
\sum_{\sigma\in Spin(M,\mathbb{Z}/\beta)}\tilde\tau_{N,K}(M,\sigma)$$
We will prove the reduction theorem for the
spin invariants.
 We need to specify some notation.
 We will use $\langle\ \rangle$ for the Homflypt
 invariant evaluated at $a$, and
$\langle\ \rangle\tilde{\ }$ for the Homflypt
 invariant evaluated at $\tilde a$.
  Note that the framing parameter does not appear
  in the coefficients of
  $\dot\Omega$,
so  that we can use them for cabling
and then evaluate the Homflypt invariant
at $a$ or $\tilde a$.

We decompose the skein element $\dot\Omega$
according to the $N\alpha$-grading (resp. the $\beta$-grading).
$$\dot\Omega=\sum_{\xi=0}^{\alpha  N-1} \dot\Omega^{(\xi)}\hspace{1cm}
\left(\text{resp. }\dot\Omega=\sum_{\nu=0}^{\beta-1} \dot\Omega_\nu\right) $$
Formulas for the modulo $\beta$ spin invariants
are then
$$\tau^{SU(N,K)}(\mathbf{S}^3(L),\sigma_c)=
\Delta^{-\sigma(L)} (\alpha^{-1}\eta)^m\langle L(\dot\Omega_{c_1},
\dots,\dot\Omega_{c_m})\rangle
$$
$$\tilde\tau_{N,K}(\mathbf{S}^3(L),\sigma_c)=
\delta^{-\sigma(L)} (\alpha^{-1}{N'}^{-1}\tilde\eta)^m
\langle L(\dot\Omega_{c_1},
\dots,\dot\Omega_{c_m})\rangle{\tilde{\ }}
$$
Here $c=(c_1,\dots,c_m)$ satisfies the modulo $\beta$
characteristic condition, and $\sigma_c$ is the modulo
$\beta$ spin structure corresponding to $c$.
We can choose $L$ so that the linking matrix is even.
(This can be considered to be a consequence of the nullity
of the cobordism group $\Omega_3^{spin}$.)
 In this case $c$ is {\em modulo $\beta$ characteristic} if
 and only if $c$ is in the kernel of the linking matrix
 modulo $\beta$.
 
\begin{lem}
We can find $\xi_1,\dots,\xi_m\in\mathbb{Z}/{\alpha N}$ such that
$\xi_i\equiv c_i\text{ mod }\beta$, for $i=1,\dots,m$
and $\xi=(\xi_1,\dots,\xi_m)$ represents an element
in the kernel of the linking matrix modulo $\alpha N$.
\end{lem}
\begin{proof}
Denote by $B$ the linking matrix.
 The Bockstein operator 
 $$Ker(B\otimes\mathbb{Z}/\beta)\rightarrow
 coker(B\otimes\mathbb{Z}/{\alpha^2N'})$$
 associated with the exact sequence
 $$0\rightarrow \mathbb{Z}/{\alpha^2N'}
 \rightarrow \mathbb{Z}/{\alpha N}
 \rightarrow \mathbb{Z}/{\beta}
 \rightarrow 0$$
 is zero since $gcd(\beta,\alpha^2N')=1$.
 \end{proof}
In the formula for $\tau^{SU(N,K)}(\mathbf{S}^3(L),\sigma_c)$
and  $\tilde\tau_{N,K}(\mathbf{S}^3(L),\sigma_c)$
we can replace $\dot\Omega_{c_i}$ by
$\sum_{j=0}^{\alpha^2 N'-1}\hat y_{K^\beta}^j \dot\Omega^{(\xi_i)}$
(with $\xi$ fixed as in the lemma above).
Using the braiding and framing coefficients for
$y_{K^\beta}$, we obtain the formula
$${\langle L(\dot\Omega_{c_1},
\dots,\dot\Omega_{c_m})\rangle}=
\sum_{j\in(\mathbb{Z}/{{\alpha^2} N'})^m}
\zeta^{{}^tjB_L j} 
\langle L(\dot\Omega^{(\xi_1)},
\dots,\dot\Omega^{(\xi_m)})\rangle$$
$${\langle L(\dot\Omega_{c_1},
\dots,\dot\Omega_{c_m})\rangle}=\alpha^{2m}
\sum_{j\in(\mathbb{Z}/{N'})^m}
\zeta^{{}^tjB_L j} 
\langle L(\dot\Omega^{(\xi_1)},
\dots,\dot\Omega^{(\xi_m)})\rangle$$
We also have
$${\langle L(\dot\Omega_{c_1},
\dots,\dot\Omega_{c_m})\rangle}{\tilde{\ }}=(\alpha^{2}N')^{m}
\langle L(\dot\Omega^{(\xi_1)},
\dots,\dot\Omega^{(\xi_m)})\rangle{\tilde{\ }}$$
For any framed link $K$,
we have $\langle K\rangle=\left(\frac{a}{\tilde a}\right)^{K.K} \langle
K\rangle\tilde{\ }$. Using that $(\frac{a}{\tilde a})^{\alpha N}=1$,
we see that  
$$\langle L(\dot\Omega^{(\xi_1)},
\dots,\dot\Omega^{(\xi_m)})\rangle=
\langle L(\dot\Omega^{(\xi_1)},
\dots,\dot\Omega^{(\xi_m)})\rangle\tilde{\ }$$
We deduce the required formula
with the following normalization of the $U(1)$ invariant
at the $N'$-th root of unity $\zeta$.
$$\tau^{U(1)}(\mathbf{S}^3(L),\zeta)=
\left(\frac{\Delta}{\delta}\right)^{-\sigma(L)}
\left(\frac{\eta}{\tilde \eta}\right)^m
\sum_{j\in(\mathbb{Z}/{N'})^m}
\zeta^{{}^tjB_L j} $$
Note that $g=\frac{\Delta}{\delta}
\frac{\eta}{\tilde \eta}$ is a Gauss sum,
whose square modulus is $g\overline g=\left(\frac{\eta}
{\tilde \eta}\right)^2=
N'$.

The other cases are obtained similarly.
 When $N+K$ is odd it is useful to note first
that in the defining expression for
$\tau^{SU(N,K)}$ we can evaluate the Homflypt invariant
at $(-a,-s,-v)$ as well as at $(a,s,v)$.
\bibliographystyle{amsplain}

\end{document}